# ON THE WEIL-ÉTALE COHOMOLOGY OF THE RING OF $S$-INTEGERS

YI-CHIH CHIU

Abstract. In this article, we first briefly introduce the history of the Weil-étale cohomology theory of arithmetic schemes and review some important results established by Lichtenbaum, Flach and Morin. Next we generalize the Weil-etale cohomology to $S$-integers and compute the cohomology for constant sheaves $\mathbb{Z}$ or $\mathbb{R}$. We also define a Weil-étale cohomology with compact support $H_c(Y_W, -)$ for $Y = \operatorname{Spec} \mathcal{O}_{F,S}$ where $F$ is a number field, and computed them. We verify that these cohomology groups satisfy the axioms state by Lichtenbaum. As an application, we derive a canonical representation of Tate sequence from $\mathrm{R}\Gamma_c(Y_W, \mathbb{Z})$. Motivated by this result, in the final part, we define an étale complex $R\,\mathbb{G}_{\mathrm{m}}$, such that the complexes $\mathrm{RHom}_{\mathbb{Z}}(\mathrm{R}\Gamma(U_{\acute{e}t}, R\,\mathbb{G}_{\mathrm{m}}), \mathbb{Z})[-2]$ and $\tau^{\leq 3}\,\mathrm{R}\Gamma_c(U_W, \mathbb{Z})$ are canonically quasi-isomorphic for arbitrary étale $U$ over $\operatorname{Spec} \mathcal{O}_F$. This quasi-isomorphism provides a possible approach to define the Weil-etale cohomology for higher dimensional arithmetic schemes, as the Weil groups are not involved in the definition of $R\,\mathbb{G}_{\mathrm{m}}$.

Contents

## 1. Introduction

Stephen Lichtenbaum conjectured that there exists a new cohomology theory, called Weil-étale cohomology, for arbitrary arithmetic schemes $X$, such that its cohomologies are connected to the zeta function of $X$. In [8], he defined a prototype of such a cohomology theory for number rings. The main idea is replacing the role of Galois groups in the étale cohomology by Weil groups. There is still no exact definition for the Weil-étale cohomology. Let $\phi : X \hookrightarrow \overline{X}$ be a fixed Nagata and Artin-Verdier style compactification, and here we list the axioms suggested by Lichtenbaum:

a) *The Weil-étale cohomology groups with compact support* $H^q_c(X,\mathbb{Z}) := H^q(\overline{X}, \phi_!\mathbb{Z})$ *are finitely generated abelian groups that equal to 0 almost everywhere, and independent of the choice of compactification of* $X$.
b) *If* $\widetilde{\mathbb{R}}$ *denotes the sheaf of real valued functions on* $X$*, then* $H^q_c(X,\widetilde{\mathbb{R}})$ *and* $H^q_c(X,\mathbb{Z})$ *are independent of the choice of* $\phi$*. Moreover, the natural map from* $H^q_c(X,\mathbb{Z})\otimes_{\mathbb{Z}}\mathbb{R}$ *to* $H^q_c(X,\widetilde{\mathbb{R}})$ *is an isomorphism.*
c) *There exists an element* $\psi \in H^1(\overline{X},\mathbb{R})$ *such that the complex*

$$\xrightarrow{\smile\psi} H^i_c(X,\widetilde{\mathbb{R}}) \xrightarrow{\smile\psi} H^{i+1}_c(X,\widetilde{\mathbb{R}}) \xrightarrow{\smile\psi}$$

*is exact. In particular,*

$$\sum_{i\geq 0}(-1)^i \dim H^i_c(X,\widetilde{\mathbb{R}}) = 0.$$

d) $ord_{s=0}\zeta_X(s) = \sum_{i\geq 0}(-1)^i i \; rank_{\mathbb{Z}}\, H^i_c(X,\mathbb{Z}).$
e) *The Euler characteristic* $\chi_c(X)$ *(see Def. 4.1) of the complex* $H^q_c(X,\mathbb{Z})$ *is well-defined, and the leading coefficient of the Taylor's expansion of* $\zeta_X(z)$ *at* $z=0$ *is* $\pm\chi_c(X)$.

The cohomology theory $H^*(\overline{X}_W,-)$ for $X = \operatorname{Spec}\mathcal{O}_F$, where $F$ is a number field, defined in [8], is an inductive limit of cohomology for a projective system of sites. And the above axioms hold only under the assumption $H^i(\overline{X}_W,\mathbb{Z}) = 0$ for all $j > 3$. However, Matthias Flach [4] proved that the cohomology group $H^i(\overline{X}_W,\mathbb{Z})$ is an infinitely generated abelian group for even $i \geq 4$, and is 0 for odd $i \geq 5$. In a recent paper [5], by modifying the prototype defined by Lichtenbaum, Flach and Baptiste Morin constructed a topos that recovers the Lichtenbaum Weil-étale cohomology groups for number fields. By this result, one may avoid the limit process for the case of number fields.

Recall that Lichtenbaum only defined the Weil-étale cohomology for $\operatorname{Spec}\mathcal{O}_F$, but one can generalize it to its arbitrary open subschemes $Y = \operatorname{Spec}\mathcal{O}_{F,S}$ in an obvious way. In section 3, we compute the Weil-étale cohomology of $Y$ with $\mathbb{Z}$-coefficients and give a brief review of Morin's results concerning Weil-étale cohomology. From section 4, we assume that the ground field $F$ is *totally imaginary.* Under this assumption, in section 4.1, we show that the Weil-étale cohomology groups with compact support for $S$-integers with $\mathbb{Z}$ or $\widetilde{\mathbb{R}}$-coefficients are as follows,

$$H^p_c(Y_W,\mathbb{Z}) = \begin{cases} 0 & p=0,\\ \prod_S \mathbb{Z}/\mathbb{Z} & p=1,\\ \mathrm{Pic}(Y)^{\mathcal{D}} \times \mathrm{Hom}(\mathrm{U_S},\mathbb{Z}) & p=2,\\ \mu_F^{\mathcal{D}} & p=3,\end{cases}$$

and

$$H^p_c(Y_W,\widetilde{\mathbb{R}}) = \begin{cases} \prod_S \mathbb{R}/\mathbb{R} & p=1,2,\\ 0 & \textit{otherwise}.\end{cases}$$

The computation is based on Theorem 4.1 that there exists an exact triangle

$$\mathrm{R}\Gamma_c(Y_{ét},\mathbb{Z}) \to \tau^{\leq 3}\,\mathrm{R}\Gamma_c(Y_W,\mathbb{Z}) \to \mathrm{Hom}(\mathrm{U_S},\mathbb{Q})[-2] \to .$$

In section 4.2, we show that all the axioms (a)-(d) hold for $S$-integers if we truncate the cohomologies of higher degrees. The proof relies on that in Lichtenbaum's original paper and a closer look of the effect of the cup product by $\psi$.

Recall that for any Galois extension $F/L$ of number fields with Galois group $G$, and a $G$-stable subset $S$ of valuation of $F$ containing all the infinite places and those ramified in $F/\mathbb{Q}$ and so that $\mathrm{Pic}(\mathcal{O}_{F,S})=0$, then, by using class field theory, one can define an exact sequence ([14] Théorèm 5.1)

$$0 \to \mathrm{U_S} \to A \to B \to X_S \to 0,$$

where $U_S$ is the group of $S$-units, $X_S := \{(x_v) \in \coprod_{v\in S}\mathbb{Z} \mid \sum x_v = 0\}$ and $A$ and $B$ are finitely generated cohomologically trivial $\mathbb{Z}[G]$-modules. We call the complex $A \to B$ the Tate sequence associated to $S$. The Tate sequence is not unique, but one can choose the complex $\Psi_S : A \to B$ so that it represents a certain canonical class of $\mathrm{Ext}^2_G(X_S,U_S)$. We denote by $\widetilde{\Psi}_S$ the complex $A\to B\to X_S\otimes_{\mathbb{Z}}\mathbb{Q}$ and $(-)^D := \mathrm{RHom}_{\mathbb{Z}}(-,\mathbb{Q}/\mathbb{Z}) = \mathrm{Hom}_{\mathbb{Z}}(-,\mathbb{Q}/\mathbb{Z})$. Burns and Flach([3] Prop. 3.1) showed that there is a map of complexes,

$$\widetilde{\Psi}_S \to R\Gamma_c(Y_{ét},\mathbb{Z})^D[-3],$$

inducing an isomorphism on $H^i$ for $i \neq 0$ and the inclusion $U_S \hookrightarrow \widehat{U_S} := U_S \otimes_{\mathbb{Z}} \widehat{\mathbb{Z}}$ on $H^0$. In section 4.3, we are able to show that $\tau^{\leq 3}R\Gamma_c(Y_W,\mathbb{Z})$ gives a canonical description of the complex $\Psi_S$ in the sense that there is an quasi-isomorphism,

$$\Psi_S \xrightarrow{qis} \mathrm{RHom}\left(\tau^{\leq 3}\,\mathrm{R}\Gamma_c\left(Y_W,\mathbb{Z}\right),\mathbb{Z}\right)[-2]. \tag{1}$$

In the final section, we construct a canonical complex $R\,\mathbb{G}_\mathrm{m}$ of $\overline{Y}_{ét}$-sheaves as an element of $\mathrm{Ext}^2_{\overline{Y}_{ét}}(\widetilde{\mathscr{X}},\mathbb{G}_\mathrm{m})$, where $\widetilde{\mathscr{X}}$ is the $\overline{Y}_{ét}$-sheaf associated to the presheaf $U_{ét} \mapsto X_{\overline{U}-U}$. We show that

$$H^i(U_{ét},R\,\mathbb{G}_\mathrm{m}) = H^i\Big(\mathrm{RHom}_{\mathbb{Z}}(\tau^{\leq 3}\,\mathrm{R}\Gamma_c(U_W,\mathbb{Z})\,,\mathbb{Z})[-2]\Big),\quad \forall i \geq 0,$$

for any étale $U \to \mathrm{Spec}\,\mathcal{O}_F$. Therefore, it may be possible to define certain Weil-etale cohomology for arbitrary arithmetic scheme by generalizing $R\,\mathbb{G}_\mathrm{m}$.

Moreover, we prove that, in fact, there is a canonical quasi-isomorphism:

$$\mathrm{R}\Gamma(U_{ét},R\,\mathbb{G}_\mathrm{m}) \xrightarrow{\mathrm{qis}} \mathrm{RHom}_{\mathbb{Z}}(\tau^{\leq 3}\,\mathrm{R}\Gamma_c(U_W,\mathbb{Z})\,,\mathbb{Z})[-2] \ \text{ in } D(\mathbb{Z}[G]),$$

for any étale $U \to \mathrm{Spec}\,\mathcal{O}_F$, when $S$ is stable under the action of $G$. When $U$ small enough so that the Tate sequence exists, this together with (1) show that

$$\mathrm{R}\Gamma(U_{ét},R\,\mathbb{G}_\mathrm{m}) \simeq \Psi_S.$$

When $G$ is the trivial group, this implies the quasi-isomorphism,

$$\mathrm{RHom}(\tau^{\leq 3}\,\mathrm{R}\Gamma_c(U_W,\mathbb{Z}),\mathbb{Z})[-2] \simeq \mathrm{R}\Gamma(U_{\acute{e}t}, R\,\mathbb{G}_{\mathrm{m}}),$$

in $D(\mathbb{Z})$, for any étale $U \to \operatorname{Spec}\mathcal{O}_F$.

Therefore, one sees that the the complex $\mathrm{R}\Gamma(U_{\acute{e}t}, R\,\mathbb{G}_{\mathrm{m}})$ generalizes the Tate sequences and its $\mathbb{Z}$-dual recovers the *truncated* Weil-etale cohomology groups of $S$-integers.

Also, this suggests the following duality property:

**Conjecture.** *There is a (perfect) paring*

$$\mathrm{RHom}_{U_{\acute{e}t}}(\mathscr{F}, R\,\mathbb{G}_{\mathrm{m}}) \times \tau^{\leq 3}\,\mathrm{R}\Gamma_c(U_W,\mathscr{F}) \to \mathbb{Z},$$

*for a certain class of étale sheaf $\mathscr{F}$ on $U$.*

The reason that it is hard to generalize Lichtenbaum's prototype to higher dimensional arithmetic schemes $\mathscr{X}$ is that there are no Weil groups for higher dimensional fields. However, Theorem 5.3 shows us a probability to generalize Lichtenbaum's prototype, because we do not use Weil groups when defining the complex $R\,\mathbb{G}_{\mathrm{m}}$. One direct conjecture is that

**Conjecture.** *For an arithmetic scheme $\mathscr{X}$ of dimension $n$, one may define a complex $R\mathbb{Z}(n)$ in $\mathrm{Ext}^2_{\overline{\mathscr{X}}}(\mathbb{Z}(n),\mathscr{F}^\bullet(n))$, so that the $\mathbb{Z}$-dual of $\mathrm{R}\Gamma(\overline{\mathscr{X}}_{et}, R\mathbb{Z}(n))$ defines certain Weil-étale cohomology theory, where $\mathscr{F}^\bullet(n)$ is a complex of étale sheaves that depends on $n$.*

## 2. Preliminary

**Notation.** Here we list the notations that are used frequently throughout this article.

- $F$ is a number field if not specified, $G_F$ is its absolute Galois group. $G$ is reserved to be the Galois group of any Galois extension $F/L$.
- $\mathcal{O}_F$: the ring of integers of $F$.
- $\overline{\operatorname{Spec}\mathcal{O}_F}$: *the Artin-Verdier compactification* of $\operatorname{Spec}\mathcal{O}_F$, which, as a set, is the same as the set of the equivalence classes of valuation of $F$ (the trivial valuation is included). We say a subset of $\overline{\operatorname{Spec}\mathcal{O}_F}$ is *open* if it has finite complement.
- For any open subset $U$ of $\overline{\operatorname{Spec}O_F}$, we denote by
  - $\overline{U} := \overline{\operatorname{Spec}O_F}$,
  - $K(U)$ : the function field of $U$,
  - $U_\infty$ (resp. $U_f$): the set of points of $U$ corresponding to infinite (resp. finite) places of $K(U)$.
- For any field $K$, $\mathcal{M}_K$ is the collection of equivalence classes of valuations of $K$.
- $v_0$: the trivial valuation of $F$.
- $S$: a finite subset of $\mathcal{M}_K$. We always require that it contain all the infinite places then $S = \overline{\operatorname{Spec}\mathcal{O}_F} - \operatorname{Spec}\mathcal{O}_{F,S}$ as sets.
- $S_\infty$(resp. $S_f$): the subset of all the infinite (resp. finite) places in $S$.
- $\mathrm{U_S}$: the group of units of $\mathcal{O}_{F,S}$.
- $(-)^\vee := \mathrm{Hom}_{\mathbb{Z}}(-,\mathbb{Z})$, $(-)^D := \mathrm{RHom}_{\mathbb{Z}}(-,\mathbb{Q}/\mathbb{Z}) = \mathrm{Hom}_{\mathbb{Z}}(-,\mathbb{Q}/\mathbb{Z})$ and $(-)^{\mathcal{D}} := \mathrm{Hom}_{cont}(-,\mathbb{R}/\mathbb{Z})$ the Pontryagin dual.

- $X_S$: $= \{(x_v) \in \coprod_{v \in S} \mathbb{Z} | \Sigma_v x_v = 0\}$. Hence $X_S^\vee = \coprod_S \mathbb{Z}/\mathbb{Z}$ and $(X_S^\vee)^\vee = X_S$.
- $C_F$: the usual idèle class group $C_F$ (resp. $F^\times$) when the field $F$ is global (resp. local).
- For any site $\mathcal{C}$, the category of abelian sheaves $Sh(\mathcal{C})$ is abelian and has enough injectives.

## 2.1. **The Weil Groups for Number Fields and Local Fields.**

2.1.1. *Weil Groups.* Idèle class groups and Galois groups are the central objects in class field theory. In 1951, Weil introduced the Weil groups, which carries the information of both. There are many equivalent ways to define the Weil groups. (cf. [13] sec.1). A brief definition of the Weil group for $F$ (global or local) is as follows : let $K/F$ be any finite Galois extension, and *the relative Weil group* $W_{K/F}$ be the extension

$$1 \to C_K \to W_{K/F} \to \mathrm{Gal}(K/F) \to 1$$

representing the canonical generator $\alpha_{K/F}$ of the group

$$H^2(\mathrm{Gal}(K/F), C_K) \cong \frac{1}{[K:F]}\mathbb{Z}/\mathbb{Z}.$$

By abuse of language, *the absolute Weil group* for $F$ is defined as

$$W_F := \varprojlim_{K/F} W_{K/F},$$

where the index runs over all finite Galois extension $K/F$. This induces a continuous morphism $\alpha_F : W_F \to G_F$ with a dense image, and we know that there is an isomorphism of topological groups $r_F : C_F \xrightarrow{\sim} W_F^{ab}$ ([1] p.238-239).

From now on we restrict $F$ to be a number field. Let $K/F$ be a finite Galois extension, and $S$ be a finite set of places of $K$ containing all the infinite ones and those ramified in $K/F$, we define $W_{K/F,S}$ to be the extension $C_{K,S} := C_K/U_{K,S}$ by $\mathrm{Gal}(K/F)$ that represents the canonical class $\alpha_{K/F,S}$ of the group

$$H^2(\mathrm{Gal}(K/F), C_{K,S}) \cong H^2(\mathrm{Gal}(K/F), C_K),$$

where $U_{K,S}$ is the idèles of the form $(a_v)$ in which $a_v = 1$ (resp. $a_v \in \mathcal{O}_{K_v}^\times$) for all $v \in S$ (resp. $v \notin S$). It is easy to see that $W_{K/F,S} \cong W_{K/F}/U_{K,S}$ and then $W_{K/F,S}^{ab} \cong C_{K,S}$. Together with the map $\log|\cdot| : C_{K,S} \to \mathbb{R}$, one defines a map $l_{L,S} : W_{K/F,S} \to \mathbb{R}$. Further, we know that $W_F = \varprojlim_{K/F,S} W_{K/F,S}$. ([8] Lemma 3.1), and there is a canonical map $l_{v_0} : W_F \to \mathbb{R}$.

2.1.2. *Weil Maps.* For any non-trivial place $v$ of $F$, we choose a place $\overline{v}$ of $\overline{F}$ lying over $v$ once for all, and denote by $D_v$ the associated decomposition group and $I_v$ the inertia group. This induces an embedding $\mathfrak{o}_v : D_v = G_{F_v} \to G_F$. We set $G_{k(v)} := D_v/I_v$. Note that then $G_{k(v)} = 1$ for any infinite place $v$. By the construction of Weil groups, there exists a so-called Weil map $\Theta_v : W_{F_v} \to W_F$, which is also an embedding such that the following diagram commutes.

$$\begin{array}{ccc} W_{F_v} & \xrightarrow{\Theta_v} & W_F \\ \downarrow{\scriptstyle \alpha_{F_v}} & & \downarrow{\scriptstyle \alpha_F} \\ G_{F_v} & \xrightarrow{\mathfrak{o}_v} & G_F. \end{array}$$

For any non-trivial place $v$, the composition $l_v : W_{F_v} \xrightarrow{p} W_{F_v}^{ab} \xrightarrow{\sim} F_v^\times \xrightarrow{|.|} \mathbb{R}$ has kernel $p^{-1}(\mathcal{O}_{F_v}^\times)$. Thus, $W_{k(v)} := W_{F_v}/p^{-1}(\mathcal{O}_{F_v}^\times) \cong \mathbb{Z}$ (resp. $\cong \mathbb{R}$) when $v \nmid \infty$ (resp. $v|\infty$.). We denote the quotient $W_{F_v} \to W_{k(v)}$ by $\pi_v$. Summarily, the diagram

$$
\begin{array}{ccc}
W_{F_v} & \xrightarrow{q_v} & W_{k(v)} \\
\downarrow{\Theta_v} & & \downarrow{l_v} \\
W_F & \xrightarrow{l_{v_0}} & \mathbb{R},
\end{array}
\tag{2}
$$

is commutative. Note that, for any finite place $v$, $l_v$ sends the generator $\sigma_v$ of $W_{k(v)}$ to $\log N(v)$.

Consider the commutative diagram

$$
\begin{array}{ccc}
W_{F_v} & \hookrightarrow & W_F \\
\downarrow & & \downarrow \\
W_{K_w/F_v} & \longrightarrow & W_{K/F,S} = W_{K/F}/U_{K,S}.
\end{array}
$$

It is not hard to show that the image of $W_{F_v}$ in $W_{K/F,S}$ is isomorphic to $W_{K_w/F_v}$ ( resp. $W_{K_w/F_v}/\mathcal{O}_{K_w}^\times$) when $v \in S$ (resp. $v \notin S$). This implies the canonical map $\pi_v : W_{F_v} \to W_{k(v)}$ factors through $\widetilde{W_{F_v}}$, where $\widetilde{W_{F_v}}$ is the image of $W_{F_v}$ in $W_{K/F,S}$. We denote $q_v : \widetilde{W_{F_v}} \to W_{k(v)}$ and $\theta_v : \widetilde{W_{F_v}} \to W_{K/F,S}$.

2.2. **The Classifying Topos and Cohomology of a Topological Group.** Let $G$ be a topological group. We denote by $B^{sm}G$ *the small classifying topos* of $G$ the category of discrete sets on which $G$ acts continuously, and by $B_{Top}G$ *the classifying site* of $G$, the category of $G$-topological spaces endowed with the local section topology $\mathcal{J}_{ls}$ (see [8] sec. 1), and $B_G := Sh(BG, \mathcal{J}_{ls})$ the topos of sheaves on this site. Let $e_G$ be the final object of $B_G$, we define the global section functor $\Gamma_G$ as $\mathrm{Hom}_{B_G}(e_G, -)$. For any abelian object $\mathcal{A}$ of $B_G$, we define the cohomology of $G$ with coefficient $\mathcal{A}$ by

$$H^i(G, \mathcal{A}) := R^i(\Gamma_G)(\mathcal{A}).$$

Let $A$ be any continuous $G$-module, and $\mathcal{A} : X \to \mathrm{Hom}_{BG}(X, A)$, then $H^i(G, \mathcal{A})$ coincides with the usual group cohomology group $H^i(G, A)$ when $G$ is a profinite group ([8] corollary 2.4).

**Theorem 2.1** ( [8] Thm. 3.6 & [4] Thm. 10.1)**.**

$$
H^q(W_F, \mathbb{Z}) = \begin{cases} \mathbb{Z} & q = 0, \\ C_F^{1,\mathcal{D}} & q = 2, \\ 0 & \textit{odd } q, \\ \textit{is of infinite rank} & \textit{even } q \geq 4, \end{cases}
$$

*where* $C_F^1$ *is the idèle class group with norm one and* $\mathcal{D} := \mathrm{Hom}_{conti}(-, \mathbb{R}/\mathbb{Z})$ *the Pontryagin dual.*

2.3. **The Étale Cohomology with Compact Support.** There are many kinds of étale cohomologies with compact support. Let $\phi : U = \operatorname{Spec}\mathcal{O}_{L,S} \to \overline{\operatorname{Spec}\mathcal{O}_L}$. We denote by $H^i_c(U_{ét}, \mathscr{F}) := H^i((\overline{\operatorname{Spec}\mathcal{O}_L})_{ét}, \phi_!\mathscr{F})$ *the cohomology with compact support*, which is equivalent to be defined by the exact triangle

$$\mathrm{R}\Gamma_c(U_{ét}, \mathscr{F}) \to \mathrm{R}\Gamma(U_{ét}, \mathscr{F}) \to \prod_{w \in S} \mathrm{R}\Gamma(L_w, \mathscr{F}) \to .$$

One can define a modified cohomology with compact support, $\widetilde{H}^i_c(U_{ét}, \mathscr{F})$, by replacing the complex $\mathrm{R}\Gamma(L_w, \mathscr{F})$ in the above exact triangle by $\mathrm{R}\Gamma_{Tate}(L_w, \mathscr{F})$ for each infinite places $w$.

Since Tate cohomology is isomorphic to the original cohomology for degree $\geq 1$, these two cohomologies with compact support are the same for degree $\geq 2$. Moreover, the difference between these two kinds of cohomologies is measured by the exact triangle,

$$\mathrm{R}\Gamma_c(U_{ét}, \mathscr{F}) \to \widetilde{\mathrm{R}\Gamma}_c(U_{ét}, \mathscr{F}) \to \bigoplus_{w \in S_\infty} \mathrm{R}\Gamma_\Delta(L_w, \mathscr{F}) \to,$$

where $\mathrm{R}\Gamma_\Delta(L_w, \mathscr{F}) := \operatorname{cone}(\mathrm{R}\Gamma(L_w, \mathscr{F}) \to \mathrm{R}\Gamma_{Tate}(L_w, \mathscr{F}))[-1]$.

Assume that $L$ is totally imaginary and set $\mathscr{F} = \mathbb{Z}$, we have the following exact triangle,

$$\mathrm{R}\Gamma_c(U_{ét}, \mathbb{Z}) \to \widetilde{\mathrm{R}\Gamma}_c(U_{ét}, \mathbb{Z}) \to \bigoplus_{w|\infty} \mathbb{Z}\,[0] \to, \tag{3}$$

2.3.1. *Artin-Verdier Duality.* Consider the pairing

$$\widetilde{\mathrm{R}\Gamma}_c(U_{ét}, \mathscr{F}) \times \mathrm{RHom}_{U_{ét}}(\mathscr{F}, \mathbb{G}_\mathrm{m}) \to \mathrm{R}\Gamma_c(U_{ét}, \mathbb{G}_\mathrm{m}) \to \mathbb{Q}/\mathbb{Z}[-3].$$

It induces a morphism

$$\widetilde{AV}(\mathscr{F}) : \mathrm{RHom}_{U_{ét}}(\mathscr{F}, \mathbb{G}_\mathrm{m}) \to \widetilde{\mathrm{R}\Gamma}_c(U_{ét}, \mathscr{F})^D[-3].$$

The Artin-Verdier duality theorem (cf. [9] Theorem 3.1) shows that $H^i(\widetilde{AV}(\mathbb{Z}))$ is an isomorphism for $i \geq 1$ and $H^0(\widetilde{AV}(\mathbb{Z})) : U_S \hookrightarrow \widehat{U}_S$, which can be summarized in an exact triangle

$$\mathrm{R}\Gamma(U_{ét}, \mathbb{G}_\mathrm{m}) \to \widetilde{\mathrm{R}\Gamma}_c(U_{ét}, \mathbb{Z})^D[-3] \to \widehat{U}_S/U_S[0] \to . \tag{4}$$

Together with the $\mathbb{Q}/\mathbb{Z}$-dual of (3), we obtain the following exact triangle

$$\mathrm{R}\Gamma(U_{ét}, \mathbb{G}_\mathrm{m}) \to \mathrm{R}\Gamma_c(U_{ét}, \mathbb{Z})^D[-3] \to \widehat{U}_S/U_S[0] \oplus (\Pi_{v|\infty}\mathbb{Q}/\mathbb{Z})[-2] \to . \tag{5}$$

We will need exact triangles (3), (4), and (5) in the last part of this article.

## 3. The Weil-Étale Cohomology of $S$-integers

3.1. **The Definition.** We fix a number field $F$ and set $Y = \operatorname{Spec}\mathcal{O}_F$. Let $K/F$ be a finite Galois extension and $S$ be a finite set of non-trivial valuations of $F$, containing all those which ramify in $K/F$. The site $T_{K/F,S}$ was defined in [8] as follows:

The objects of $T_{K/F,S}$ are the collections

$$((X_v), (f_v))_{v \in \overline{Y}},$$

where $X_v$ is an $W_{k(v)}$-space, and $f_v : X_v \to X_{v_0}$ is map of $W_{F_v}$-spaces. ( We regard $X_v$ as an $W_{F_v}$-space via $\pi_v$, and $X_{v_0}$ as an $W_{F_v}$-space via $\Theta_v$). Further, we require that the action of $W_F$ on $X_{v_0}$ factors through $W_{L/F,S}$.

Let $\mathcal{X} = ((X_v), (f_v))$ and $\mathcal{X}' = ((X'_v), (f'_v))$ be objects of $T_{K/F,S}$, then $Hom(\mathcal{X}, \mathcal{X}')$ is the collections of $W_{F_v}$-maps $g_v : X_v \to X'_v$ such that $g_{v_0} f_v = f'_{v_0} g_v$ for all $v$. And the fibre product of two morphisms with the same codomain is defined componentwise.

The covering $Cov(T_{K/F,S})$ consists of the family of morphisms $\{((X_{i,v}), (f_{i,v})) \to ((X_v), (f_v))\}_i$ such that $\{X_{i,v} \to X_v\}_i$ is a local section covering of $X_v$ for all $v$. We denote $\widetilde{T_{K/F,S}}$ the topos $Sh(T_{K/F,S}, \mathcal{J}_{ls})$.

Clearly, $*_{K/F,S} = ((pt), (id))$ is the final object of $T_{K/F,S}$. For any abelian sheaf $\mathscr{F}_{K,S}$ on it, we define the cohomology $H^p(T_{K/F,S}, \mathscr{F}_{K,S})$ by $H^p(T_{K/F,S}, *, \mathscr{F}_{K,S})$.

For each non-trivial place $v$, the map $i^{-1}_{K,S,v} : ((X_v)) \mapsto X_v$ induces an embedding of topoi $i_{K,S,v} : B_{W_{k(v)}} \to \widetilde{T_{K/F,S}}$. When $v = v_0$, we denoted by $j_{K,S}$ the embedding $i_{K,S,v} : B_{W_{K/F,S}} \to \widetilde{T_{K/F,S}}$.

For any finite Galois extension $K'$ of $K$, and $S' \supset S$, the canonical morphism $p : W_{K'/F,S'} \to W_{K/F,S}$ defines a transition map $t_{K'/K,S'/S} : \widetilde{T_{K'/F,S'}} \to \widetilde{T_{K/F,S}}$ ( by regarding any $W_{K'/F,S'}$-space as an $W_{K/F,S}$-space via $p$). We will throw out $K, S$ from the index if there is no risk of confusing.

**Proposition 3.1** ([11] Prop. 3.5, 3.14)**.** a) *$i_v^*$, $i_{v,*}$ and $j^*$ are exact, and $i_{v,*} i_v^* = id$.*

b) *We have the following commutative diagram*

$$\begin{array}{ccccc} B_{W_{K'/F,S'}} & \xrightarrow{j_{K',S'}} & \widetilde{T_{K'/F,S'}} & \xleftarrow{i_{K',S',v}} & B_{W_{k(v)}} \\ \downarrow & & \downarrow t & \swarrow i_{K,S,v} & \\ B_{W_{K/F,S}} & \xrightarrow{j_{K,S}} & \widetilde{T_{K/F,S}}. & & \end{array}$$

**Proposition 3.2** (ibid., Prop. 3.15, 3.16)**.** a) *There is a morphism*

$$f_{K,S} : \widetilde{T_{K/F,S}} \to B_{\mathbb{R}}$$

*so that $f_{K,S} \circ i_v$ is isomorphic to $B_{l_v} : B_{W_{k(v)}} \to B_{\mathbb{R}}$, for any closed point $v$ of $\overline{Y}$.*

b) *The following diagram is commutative*

$$\begin{array}{ccc} \widetilde{T_{K'/F,S'}} & \longrightarrow & \widetilde{T_{K/F,S}} \\ & \searrow & \downarrow \\ & & B_{\mathbb{R}} \end{array}$$

*for any $\overline{F}/K'/K/F$ and $S \subset S'$.*

**Definition 3.1.** *A* **compatible system of abelian sheaves on the sites** $(T_{K/F,S})_{K,S}$ *is a collection $\{(F_{K,S}), (\phi_{K'/K,S'/S})\}$, where $F_{K,S} \in \widetilde{T_{K/F,S}}$ and $\phi_{K'/K,S'/S}$ is a morphism $t^*_{K'/K,S'/S} F_{K,S} \to F_{K',S'}$, such that if $K''/K'/K$ a series of finite Galois extensions and $S'' \supset S' \supset S$, then $\phi_{K''/K,S''/S} = t^*_{K''/K',S''/S'}(\phi_{K''/K',S''/S'}) \circ \phi_{K'/K,S'/S}$.*

**Definition 3.2.** *Let $\mathscr{F} = \{(F_{K,S}), (\phi_{K'/K,S'/S})\}$ be any compatible system of sheaves on the sites $(T_{K/F,S})_{K,S}$, then we define the Weil-étale cohomology on* $\operatorname{Spec}\mathcal{O}_F$ *for $\mathscr{F}$, $H^p(\overline{Y}_W, \mathscr{F})$ to be $\varinjlim_{K,S} H^p(T_{K/F,S}, F_{K,S})$.*

**Example 3.1.** *Let $p_{K,S}$ be the morphism $B_{W_F} \to B_{W_{K/F,S}}$ and $\mathcal{A}$ be an abelian object of $B_{W_F}$ associated to a topological abelian $W_F$-group $A$. Then $\widetilde{A} := \{(j_{K,S,*}p_{K,S,*}\mathcal{A}), (\phi)\}$ forms a compatible system of sheaves, where*

$$\phi_{K'/K,S'/S} : t^* j_{K,S,*} p_{K,*}\mathcal{A} = t^* t_* j_{K',S',*} p_{K',*}\mathcal{A} \to j_{K',S',*} p_{K',*}\mathcal{A}$$

*is the natural map induced by the adjunction. We define $H^p(\overline{Y}_W, A)$ to be $H^p(\overline{Y}_W, \widetilde{A})$. In particular, if $A$ is a topological abelian group on which $W_F$ acts trivially, then $\widetilde{A}$ is called the* **constant sheaf** *with $A$-value. We denote by $A$ the constant sheaf defined by $A$ if $A$ has the discrete topology.*

In [11], Morin gave a direct description of the topos $\widetilde{T_{K/F,S}}$. Recall that $\widetilde{W_{K_v}}$ denotes the image of $W_{K_v}$ in $W_{K/F,S}$, and $\theta_v : \widetilde{W_{K_v}} \to W_{K/F,S}$ and $q_v : \widetilde{W_{K_v}} \to W_{k(v)}$ are the induced continuous maps.

**Definition 3.3** (ibid., Def. 3.1)**.** *We define a category $\mathfrak{F}_{K/F,S}$ as follows. The objects of this category are of the form $\mathcal{F} = (F_v; f_v)_{v \in \overline{Y}}$, where $F_v$ is an object of $B_{W_{k(v)}}$ for $v \neq v_0$ ( resp. of $B_{W_{K/F,S}}$ for $v = v_0$) and*

$$f_v : q_v^*(F_v) \to \theta_v^*(F_{v_0})$$

*is a morphism of $B_{\widetilde{W_{K_v}}}$ so that $f_{v_0} = id_{F_{v_0}}$. A morphism $\phi$ from $\mathcal{F} = (F_v; f_v)_{v\in\overline{Y}}$ to $\mathcal{F}' = (F'_v; f'_v)_{v\in\overline{Y}}$ is a family of morphisms $\phi_v : F_v \to F'_v \in Fl(B_{W_{K_v}})$(and $\phi_{v_0} \in Fl(B_{W_{K/F,S}})$) so that*

$$\begin{array}{ccc} q_v^*(F_v) & \xrightarrow{q_v^*(\phi_v)} & q_v^*(F'_v) \\ \downarrow f_v & & \downarrow f'_v \\ q_v^*(F_{v_0}) & \xrightarrow{\theta_v^*(\phi_{v_0})} & \theta_v^*(F'_{v_0}) \end{array}$$

*is a commutative diagram of $B_{\widetilde{W_{K_v}}}$.*

*For $K = \overline{F}$ and $S$ the set of all non-trivial valuations of $F$, one has $W_{K/F,S} = W_F, \widetilde{W_{F_v}} = W_{F_v}$ and we set $\mathfrak{F}_{K/F,S} = \mathfrak{F}_{W;\overline{Y}}$.*

**Proposition 3.3** (ibid., Thm. 5.9)**.** *$\mathfrak{F}_{K/F,S}$ is a topos and is equivalent to $\widetilde{T_{K/F,S}}$ as a topos.*

In the following, we identify $\mathfrak{F}_{K/F,S}$ and $\widetilde{T_{K/F,S}}$.

**Proposition 3.4** (a special case of [11] Lemma 4.5 & Prop. 4.6)**.** *Let $\mathcal{A}$ be an abelian object in $B_{W_F}$ defined by $A$ on which $W_F$ acts trivially, then*

a) *$R^p j_*\mathcal{A} := \{(R^p(j_{K,S,*})(p_{K,S,*}\mathcal{A}), (\phi)\}$ is a compatible system of sheaves.*
b) *There exists a spectral sequence*

$$H^p(\overline{Y}_W, R^q j_*\mathcal{A}) \Longrightarrow H^{p+q}(W_F, A).$$

**Proposition 3.5** ([8] Thm 5.10, [4] Theorem 11.1)**.**

$$H^p(\overline{Y}_W, \mathbb{Z}) = \begin{cases} \mathbb{Z} & p = 0, \\ 0 & p = 1, \\ \mathrm{Pic}^1(\overline{Y})^{\mathcal{D}} & p = 2, \\ \mu_F^D & p = 3, \\ \textit{is of infinite rank} & \textit{even } p \geq 4, \\ 0 & \textit{odd } p \geq 5, \end{cases}$$

*where* $\mathrm{Pic}^1(\overline{Y})$ *is the kernel of the absolute value map from* $\mathrm{Pic}(\overline{Y})$ *to* $\mathbb{R}^{>0}$ *and* $\mathrm{Pic}(\overline{Y})$ *is the Arakelov Picard group of F, which is obtained by taking the idèle group of F and dividing by the principal idèles and the unit idèles( i.e. those idèle* $(u_v)$ *s.t.* $|u_v|_v = 1$ *for all v.).*

**Proposition 3.6** ([8] Theorem 5.11)**.**

$$H^p(\overline{Y}_W, \widetilde{\mathbb{R}}) = \begin{cases} \mathbb{R} & p = 0, 1, \\ 0 & p > 1. \end{cases}$$

More generally, for any open subscheme $U = \operatorname{Spec} \mathcal{O}_{F,S}$ of $\overline{Y}$, $U$ defines an open sub-topos $\mathfrak{F}_{K/F,S'}/_U \to \mathfrak{F}_{K/F,S'}$, the full sub-category of $\mathfrak{F}_{K/F,S'}$ whose objects are of the form $\mathcal{F} = (F_v; f_v)_{v \in U}$ (i.e. $F_v = \emptyset$ for $v \in \overline{Y} - U$). One can show that $\mathfrak{F}_{K/F,S'}/_U$ is equivalent to $T^U_{K/F,S'}$, where $T^U_{K/F,S'}$ is defined similarly as $T_{K/F,S'}$ by throwing out from each object all the $v$-components for $v \in \overline{Y} - U$. We have the usual embeddings $j^S_{K/F,S'} : B_{W_{K/F,S'}} \to \mathfrak{F}_{K/F,S'}/_U$ and $i^S_v : B_{W_{k(v)}} \to \mathfrak{F}_{K/F,S'}/_U$.

**Definition 3.4.** *We define the cohomology* $H^p(\mathfrak{F}_{K/F,S'}, U, -) := H^p(\mathfrak{F}_{K/F,S'}/_U, -)$, *and* $H^p(U_W, \mathbb{Z}) = \varinjlim H^p(\mathfrak{F}_{K/F,S'}, U, \mathbb{Z})$.

**Proposition 3.7.** *For any proper open subscheme* $U = \operatorname{Spec} \mathcal{O}_{F,S}$ *of* $\overline{Y}$*,*

$$H^p(U_W, \mathbb{Z}) = \begin{cases} \mathbb{Z} & p = 0, \\ 0 & p \textit{ is odd,} \\ (C^1_F / \coprod_{v \in U} U_v)^{\mathcal{D}} & p = 2, \\ H^p(W_F, \mathbb{Z}) & p \geq 3. \end{cases}$$

*Proof.* There exists the Leray spectral sequence induced by the inclusion of the generic point

$$\varinjlim_{K/F,S'} H^p(\mathfrak{F}_{K/F,S'}/_U, R^q j^S_{K/F,S',*}\mathbb{Z}) \Longrightarrow H^{p+q}(W_F, \mathbb{Z}).$$

By the same argument of [8] Theorem 4.8, for $q > 0$, one has $R^q j^S_{K/F,S',*}\mathbb{Z} = \coprod_{v \in S' \setminus S} i^S_{v\,*} i^{S*}_v R^q j^S_{K/F,S',*}\mathbb{Z}$. Moreover, the functor $i^{S*}_v R^q j^S_{K/F,S',*} = R^q(i^{S*}_v j^S_{K/F,S',*}) = R^q(i^{S*}_v i^{S*} i^S{}_* j^S_{K/F,S',*}) = R^q(i^*_v j_{K/F,S',*}) = i^*_v R^q j_{K/F,S',*}$, where $i^S$ is the open embedding $\mathfrak{F}_{K/F,S'}/_U \to \mathfrak{F}_{K/F,S'}$.

Therefore, for any $q > 0$, $H^p(\mathfrak{F}_{K/F,S'}/_U, R^q j^S_{K/F,S',*}\mathbb{Z}) = \Pi_{v \in S' \setminus S} H^p(W_{k(v)}, i^*_v R^q j_{K/F,S',*}\mathbb{Z})$.

By the computation due to Flach [4], we have

$$\varinjlim_{K/F,S'} H^p(\mathfrak{F}_{K/F,S'}/_U, R^q j^S_{K/F,S',*}\mathbb{Z}) = \begin{cases} \mathbb{Z} & q = 0, \\ \prod_{v \in U} U_v^{\mathcal{D}} & p = 0 \text{ and } q = 2, \\ 0 & p > 0 \text{ and } q \neq 2. \end{cases}$$

So the 2nd page of the spectral sequence looks like the following,

$$\begin{array}{cccc} 0 & & & \\ \prod_{v\in U} U_v^{\mathcal{D}} & & 0 & \\ 0 & & & \\ H^0(U_W,\mathbb{Z}) & H^1(U_W,\mathbb{Z}) & H^2(U_W,\mathbb{Z}) & H^3(U_W,\mathbb{Z}). \end{array}$$

This gives us $H^0(U_W,\mathbb{Z}) = \mathbb{Z}$, $H^1(U_W,\mathbb{Z}) = H^1(W_F,\mathbb{Z}) = 0$ and the exact sequence,

$$0 \longrightarrow H^2(U_W,\mathbb{Z}) \longrightarrow H^2(W_F,\mathbb{Z}) \longrightarrow \prod_{v\in U} U_v^{\mathcal{D}} \longrightarrow H^3(U_W,\mathbb{Z}) \longrightarrow H^3(W_F,\mathbb{Z}) = 0,$$
$$H^2(W_F,\mathbb{Z}) = C_F^{1,\mathcal{D}}$$

which is dual to

$$0 \longrightarrow H^3(U_W,\mathbb{Z})^{\mathcal{D}} \longrightarrow \coprod_{v\in U} U_v \longrightarrow C_F^1 \longrightarrow H^2(U_W,\mathbb{Z})^{\mathcal{D}} \longrightarrow 0.$$

As the middle map is injective when $U$ is a proper subset of $\overline{Y}$, we see that $H^2(U_W,\mathbb{Z}) = (C_F^1/\coprod_{v\in U} U_v)^{\mathcal{D}}$ and $H^3(U_W,\mathbb{Z}) = 0$. And for $p \geq 3, H^p(U_W,\mathbb{Z}) = H^p(W_F,\mathbb{Z})$, which is vanishing for odd $p \geq 3$. □

**Definition 3.5.** *Let $\mathscr{F} = (F_{K/F,S})$ be a compatible system on sheaves on the site $(T_{K/F,S})$. Then $i_v^*\mathscr{F} := (i_{K,S,v}^* F_{K/F,S})$ a compatible system of sheaves on the site $B_{W_{k(v)}}$ and it gives rise to restriction maps $H^*(T_{K/F,S}, F_{K/F,S}) \to H^*(k(v)_W, i_{K,S,v}^* F_{K/F,S})$. Let $U$ be any open subscheme of $\overline{Y}$. Since $\overline{U} = \overline{Y}$, we define the* **Weil-étale cohomology with compact support** *for an open subscheme $U$ with coefficient $\mathscr{F}$ as the cohomologies of the mapping cone*

$$\mathrm{R}\Gamma_c(U_W,\mathscr{F}) := \mathrm{cone}\Big(\mathrm{R}\Gamma(\overline{Y}_W,\mathscr{F}) \to \bigoplus_{v\in\overline{Y}\setminus U} \mathrm{R}\Gamma(k(v)_W, i_v^*\mathscr{F})\Big)[-1],$$

*where $k(v)_W := B_{W_{k(v)}}$. The map $\mathrm{R}\Gamma(\overline{Y}_W,\mathscr{F}) \to \mathrm{R}\Gamma(k(v)_W, i_v^*\mathscr{F})$ is derived by first applying global section functor to $I^\bullet_{F_{K/F,S}} \to i_{K,S,v,*} I^\bullet_{i^*_{K,S,v}F_{K/F,S}}$ and then taking the inductive limits, where $I^\bullet_{F_{K/F,S}}$ (resp. $I^\bullet_{i^*_{K,S,v}F_{K/F,S}}$) is an injective resolution of $F_{K/F,S}$(resp. $i^*_{K,S,v}F_{K/F,S}$). Note that the complex $\mathrm{R}\Gamma_c(U_W,\mathscr{F})$ is uniquely determined up to homotopy equivalences.*

As $H^p(k(v)_W,\mathbb{Z}) = H^p(W_{k(v)},\mathbb{Z})$ has non-vanishing cohomology (all equals to $\mathbb{Z}$) only when $v|\infty$ and $p = 0$ or $v \nmid \infty$ and $p = 0, 1$. And $H^p(k(v)_W,\mathbb{R})$ equals to $\mathbb{R}$ (resp. 0) if $p = 0, 1$ (resp. $p > 1$). One can easily deduce the following,

**Proposition 3.8** ([8] Thm. 6.3)**.** *Let $Y = \operatorname{Spec} \mathcal{O}_F$, then*

$$H_c^p(Y_W,\mathbb{Z}) = \begin{cases} 0 & p = 0, \\ (\coprod_{v|\infty}\mathbb{Z})/\mathbb{Z} & p = 1, \\ H^p(\overline{Y}_W,\mathbb{Z}) & p \geq 2, \end{cases}$$

*and*

$$H^p_c(Y_W,\mathbb{R}) = \begin{cases} (\coprod_{v|\infty}\mathbb{R})/\mathbb{R} & p=1,2, \\ 0 & \textit{otherwise.} \end{cases}$$

3.2. **The Weil-Étale Sites and the Artin-Verdier Étale Sites.** All the material in this subsection can be found in [11]. We will only give brief definitions and state the results.

**Proposition 3.9** (ibid., Prop 7.1)**.** *There exists a morphism of left exact sites*

$$\begin{array}{llll} \zeta^*: & (Et(\overline{Y});\mathcal{J}_{\acute{e}t}) & \to (T_{W;\overline{Y}};\mathcal{J}_{ls}) \\ & \quad\overline{X} & \mapsto \mathcal{X} = ((X_v);(\varphi_v)), \end{array}$$

*where* $X_v := Hom_{\overline{Y}}(\overline{Y}^{sh}_v;\overline{X})$ *and* $\overline{Y}^{sh}_v = \mathrm{Spec}(\overline{k(v)})$ *(resp.* $=(\emptyset;\overline{v})$*) when* $v\nmid\infty$ *(resp.* $v|\infty$*).*

Let $K/F$ be a finite Galois extension. We denote by $Et_{K/F}$ the full sub-category of $Et(\overline{Y})$ consisting of étale $\overline{Y}$-schemes $\overline{X}$ such that the action of $G_F$ on the finite set $X_{v_0} := Hom(\mathrm{Spec}(\overline{F});\overline{X})$ factors through $\mathrm{Gal}(K/F)$. This category is again endowed with the topology induced by the étale topology on $Et(\overline{Y})$ via the inclusion functor $Et_{K/F}\to Et(\overline{Y})$. These morphisms are compatible and induce a projective system $(\widetilde{(Et_{K/F},\mathcal{J}_{\acute{e}t})})_K$.

**Proposition 3.10** (ibid., Prop. 7.7)**.** *The canonical morphism*

$$\overline{Y}_{\acute{e}t}\to \varprojlim\, \widetilde{(Et_{K/F},\mathcal{J}_{\acute{e}t})}$$

*is an equivalence.*

**Proposition 3.11** (ibid., Prop. 7.5 & 7.8)**.** a) $\zeta^*:(Et(\overline{Y});\mathcal{J}_{\acute{e}t})\to(T_{W;\overline{Y}};\mathcal{J}_{ls})$ *is the inverse image of a morphism of topoi* $\zeta:\mathfrak{F}_{W;\overline{Y}}\to\overline{Y}_{\acute{e}t}$*.*

b) $\zeta$ *induces a morphism of topoi* $\zeta_{K,S}:\mathfrak{F}_{K/F,S}\to\widetilde{(Et_{K/F},\mathcal{J}_{\acute{e}t})}$*. Moreover, the diagram*

$$\begin{array}{ccc} \mathfrak{F}_{K'/F,S'} & \xrightarrow{\zeta_{K',S'}} & \widetilde{(Et_{K'/F},\mathcal{J}_{\acute{e}t})} \\ \downarrow & & \downarrow \\ \mathfrak{F}_{K/F,S} & \xrightarrow{\zeta_{K,S}} & \widetilde{(Et_{K/F},\mathcal{J}_{\acute{e}t})} \end{array}$$

*is commutative, for* $K'/K/F$ *and* $S\subset S'$*.*

**Proposition 3.12** (ibid., Prop. 4.69)**.** *For any closed point* $v\in\overline{Y}$*, the diagram*

$$\begin{array}{ccc} B_{W_{k(v)}} & \xrightarrow{\alpha_v} & B^{sm}_{G_{k(v)}} \\ \downarrow{\scriptstyle i_v} & & \downarrow{\scriptstyle u_v} \\ \mathfrak{F}_{W,\overline{Y}} & \xrightarrow{\zeta} & \widetilde{\overline{Y}_{\acute{e}t}}. \end{array}$$

*is a pull-back of topoi.*

**Theorem 3.1** (ibid., Thm. 8.5 & Prop. 8.6)**.** *Let $\mathscr{F} = (F_{K,S}; f_t)$ be a compatible system of sheaves on the sites $(T_{K/F,S})_{K,S}$ . There exists a bounded below complex $R\mathscr{F}$ of abelian $\overline{Y}_{ét}$-sheaves and an isomorphism*

$$H^p(\overline{Y}_{ét}, R\mathscr{F}) \simeq H^p(\overline{Y}_W, \mathscr{F}),$$

*where the left-hand side is the étale hypercohomology of the complex $R\mathscr{F}$. In particular, one has a spectral sequence relating Lichtenbaum's Weil-étale cohomology to étale cohomology,*

$$H^p(\overline{Y}_{ét}, R^q\mathscr{F}) \Longrightarrow H^{p+q}(\overline{Y}_W, \mathscr{F}),$$

*where $R^q\mathscr{F} := H^q(R\mathscr{F})$. The complex $R\mathscr{F}$ is well-defined up to quasi-isomorphism and functorial in $\mathscr{F}$. Moreover, $R^q\mathscr{F}$ is the sheaf associated to the presheaf*

$$\begin{aligned} \mathcal{P}^q\mathscr{F}: \quad \overline{Y}_{ét} &\longrightarrow \underline{Ab} \\ U &\mapsto H^q(U_W, \mathscr{F}). \end{aligned}$$

**Remark.** *For the future use, we describe the construction of the complex $R\mathcal{A}$.*

*For each pair $(K,S)$, we fix an injective resolution $0 \to F_{K,S} \to I^\bullet_{K,S}$ in $\mathfrak{F}_{K,S}$. and with morphism of complexes $t^* I^\bullet_{K,S} \to I^\bullet_{K',S'}$, for any transition map $t$. These complexes are compatible. Then one can show that $(\zeta_{K,S,*} I^\bullet_{K,S})_S$ form a direct systems of complexes. We set*

$$I^\bullet_K := \varinjlim_S \zeta_{K,S,*} I^\bullet_{K,S}.$$

*Then $(I^\bullet_K)_K$ defines a compatible system of complexes of sheaves on the sites $(Et_{K/F})_K$. Finally, the complex*

$$R\mathscr{F} := \varinjlim_K u_K^* I^\bullet_K,$$

*where $u : \overline{Y}_{ét} \to (\widetilde{Et_{K/F}}, \mathcal{J}_{ét})$.*

We are interested in the complex $\tau^{\leq 3}\,\mathrm{R}\Gamma(\overline{Y}_W, \mathbb{Z})$, and Morin proved the following.

**Theorem 3.2** (ibid., Thm 8.5 & 9.5)**.** *Suppose $F$ is totally imaginary then the following hold*

a) $\mathrm{R}\Gamma(\overline{Y}_{ét}, \tau^{\leq 3} R\mathbb{Z}) \xrightarrow{qis} \tau^{\leq 3}\,\mathrm{R}\Gamma(\overline{Y}_W, \mathbb{Z})$,
b) $\mathrm{R}\Gamma(\overline{Y}_{ét}, R^2\mathbb{Z}) \xrightarrow{qis} \mathrm{Hom}_{\mathbb{Z}}(U_F, \mathbb{Q})[0]$.

Here is an immediate corollary of this theorem:

**Corollary 3.1.** *We have the following exact triangle*

$$\mathrm{R}\Gamma(\overline{Y}_{ét}, \mathbb{Z}) \to \tau^{\leq 3}\,\mathrm{R}\Gamma(\overline{Y}_W, \mathbb{Z}) \to \mathrm{Hom}_{\mathbb{Z}}(U_F, \mathbb{Q})[-2] \to .$$

*Proof.* By Prop. 3.7, $H^1(U_W, \mathbb{Z}) = 0$ for all open étale $U$ over $\overline{Y}_{ét}$, so $R^1\mathbb{Z} = 0$ by Theorem 1.2. Thus, we have an exact triangle

$$\mathbb{Z}[0] \longrightarrow \tau^{\leq 3} R\mathbb{Z} \longrightarrow R^2\mathbb{Z}[-2] \longrightarrow .$$

Applying $\mathrm{R}\Gamma(\overline{Y}_{ét}, -)$ to the above exact triangle and by using the previous Theorem, we obtain the desired exact triangle. □

## 4. Weil-étale cohomology with compact support

In this whole section, we assume that $F$ is a *totally imaginary* number field, and $Y = \mathrm{Spec}\,\mathcal{O}_{F,S}$.

4.1. **The Computation of the Cohomology with Compact Support.** Let $A$ be an abelian group with trivial $W_F$-action and $A$ the constant sheaf defined by it. Recall that (cf. Def. 3.5) the Weil-étale cohomology groups with compact support $H_c^*(Y_W, A)$ are defined by the cohomologies of the mapping cone

$$\mathrm{R}\Gamma_c(Y_W, A) = \mathrm{cone}\Big(\mathrm{R}\Gamma(\overline{Y}_W, A) \to \bigoplus_{v \in S} \mathrm{R}\Gamma(k(v)_W, A)\Big)[-1],$$

which is equivalent to the exact triangle

$$\mathrm{R}\Gamma_c(Y_W, A) \to \mathrm{R}\Gamma(\overline{Y}_W, A) \to \bigoplus_{v \in S} \mathrm{R}\Gamma(k(v)_W, A) \to .$$

We are interested in the complex $\tau^{\leq 3}\,\mathrm{R}\Gamma_c(Y_W, \mathbb{Z})$, and want to prove first the following.

**Theorem 4.1.** *There is an exact triangle*

$$\mathrm{R}\Gamma_c(Y_{ét}, \mathbb{Z}) \longrightarrow \tau^{\leq 3}\,\mathrm{R}\Gamma_c(Y_W, \mathbb{Z}) \longrightarrow \mathrm{Hom}(\mathrm{U_S}, \mathbb{Q})[-2] \to . \tag{6}$$

Before proving Theorem 4.1, we show the following lemma.

**Lemma 4.1.** *For any constant sheaf $A$, the inclusion $A[0] \to RA$ induces a morphism $\mathrm{R}\Gamma(\overline{Y}_{ét}, A) \longrightarrow \mathrm{R}\Gamma(\overline{Y}_W, A)$, and the following diagram commutes*

$$\begin{array}{ccc} \mathrm{R}\Gamma(\overline{Y}_{ét}, A) & \longrightarrow & \mathrm{R}\Gamma(k(v)_{ét}, A) \\ \downarrow & & \downarrow{\scriptstyle \alpha_v^*} \\ \mathrm{R}\Gamma(\overline{Y}_W, A) & \longrightarrow & \mathrm{R}\Gamma(k(v)_W, A), \end{array}$$

*where the morphisms of the rows are the canonical ones and $\alpha_v^*$ is induced by the canonical morphism of topos $\alpha_v : B_{W_{k(v)}} \to B^{sm}_{G_{k(v)}}$.*

*Hence, taking the truncation functor $\tau^{\leq 3}$ on every terms gives rise to the following commutative diagram.*

$$\begin{array}{ccc} \tau^{\leq 3}\,\mathrm{R}\Gamma(\overline{Y}_{ét}, A) & \longrightarrow & \mathrm{R}\Gamma(k(v)_{ét}, A) \\ \downarrow & & \downarrow{\scriptstyle \alpha_v^*} \\ \tau^{\leq 3}\,\mathrm{R}\Gamma(\overline{Y}_W, A) & \longrightarrow & \mathrm{R}\Gamma(k(v)_W, A). \end{array}$$

*Proof.* Recall that by the knowledge of section 1, there exists a commutative diagram of topoi

$$\begin{array}{ccc} k(v)_W & \xrightarrow{\alpha_v} & k(v)_{ét} \\ \downarrow{\scriptstyle i_v} & & \downarrow{\scriptstyle i} \\ \widetilde{T_{\overline{Y}}} & \xrightarrow{\zeta} & \widetilde{Et_F} \\ \downarrow & & \downarrow{\scriptstyle u_K} \\ \widetilde{T_{K/F,S}} & \xrightarrow{\zeta_{K,S}} & \widetilde{Et_{K/F}}, \end{array} \qquad i_{K,S} : k(v)_W \to \widetilde{T_{K/F,S}}$$

where $\alpha_v : B_{W_{k(v)}} \to B^{sm}_{G_{k(v)}}$. Here we drop $v$ in the indexes.

Clearly, the following diagram commutes.

$$\begin{array}{ccc} & \mathrm{R}\Gamma(\overline{Y}_{\acute{e}t}, A) & \longrightarrow & \mathrm{R}\Gamma(\overline{Y}_{\acute{e}t}, i_*i^*A) = \mathrm{R}\Gamma(k(v)_{\acute{e}t}, A) \\ & \downarrow & & \downarrow \\ \mathrm{R}\Gamma(\overline{Y}_W, A) = & \mathrm{R}\Gamma(\overline{Y}_{\acute{e}t}, RA) & \longrightarrow & \mathrm{R}\Gamma(\overline{Y}_{\acute{e}t}, i_*i^*RA) \end{array} .$$

To complete the proof, we need to show the morphism $\mathrm{R}\Gamma(\overline{Y}_W, A) \to \mathrm{R}\Gamma(k(v)_W, A)$ factors through $\mathrm{R}\Gamma(\overline{Y}_W, A) \to \mathrm{R}\Gamma(\overline{Y}_{\acute{e}t}, i_*i^*RA)$.

Note that the morphism $\mathrm{R}\Gamma(\overline{Y}_W, A) \to \mathrm{R}\Gamma(k(v)_W, A)$ is induced by

$$\varinjlim \mathrm{R}\Gamma(\widetilde{T_{K/F,S}}, A) \longrightarrow \varinjlim \mathrm{R}\Gamma(\widetilde{T_{K/F,S}}, i_{K,S,*}i^*_{K,S}A).$$

By Theorem 3.1, this is the same as the morphism of complexes of $\overline{Y}_{\acute{e}t}$-sheaves

$$RA := \varinjlim u_K^* \varinjlim \zeta_{K,S,*} I^\bullet_{K,S} \longrightarrow R'A := \varinjlim u_K^* \varinjlim \zeta_{K,S,*} i_{K,S,*} J^\bullet_{K,S},$$

where $I^\bullet_{K,S}$ (resp. $J^\bullet_{K,S}$) is any fixed injective resolution of $A$ (resp. $i^*_{K,S}A$) and note that $i_{K,S,*}$ is exact and preserves injectives.

By applying the natural functor $id \to i_*i^*$ to $RA \to R'A$, we get the following canonical commutative diagram

$$\begin{array}{ccc} RA & \longrightarrow & R'A \\ {\scriptstyle p_{RA}}\downarrow & & \downarrow{\scriptstyle p_{R'A}} \\ i_*i^*RA & \longrightarrow & i_*i^*R'A. \end{array}$$

Note that

$$\begin{aligned} R'A &= \varinjlim_K u_K^* \varinjlim_S \zeta_{K,S,*} i_{K,S,*} J^\bullet_{K,S} \\ &= \varinjlim_K u_K^* \varinjlim_S u_{K,*} i_* \alpha_{v,*} J^\bullet_{K,S} \\ &= \varinjlim_K \varinjlim_S u_K^* u_{K,*} i_* \alpha_{v,*} J^\bullet_{K,S} \\ &= \varinjlim_K \varinjlim_S i_* \alpha_{v,*} J^\bullet_{K,S} \quad (\varinjlim_K \varinjlim_S u_K^* u_{K,*} F_{K,S} = \varinjlim_K \varinjlim_S F_{K,S} \ \ (\textit{c.f.}\ [10]\ \text{lemma 7.5})) \\ &= i_* R\alpha_{v,*} A. \end{aligned}$$

Because applying with $i^*$ induces the adjoint isomorphism $Hom(i_*P, i_*Q) = Hom(P, Q)$, and $i^*p_{R'A} = id$, by definition, so $p_{R'A} = id$ and the result follows. □

**Remark.** *Whenever there is a square commutative diagram of complexes in an abelian category $\mathcal{A}$, then by taking mapping cones, one can complete it as a semi-commutative 3 by 3 diagram of complexes in which rows and columns are exact triangles in the derived category $D(\mathcal{A})$. (c.f.* [16] *Exercise 10.2.6).*

*Proof of Theorem 4.1.* By [10] Prop. 6.5 and 6.6,

$$R^q\alpha_{v,*}\mathbb{Z} = \begin{cases} \mathbb{Z}, \mathbb{Q}, 0 & \text{for } q = 0, 1, \text{ and } q \geq 2 \ (v \nmid \infty), \\ \mathbb{Z}, 0 & \text{for } q = 0, \text{ and } q \geq 1 \ (v | \infty). \end{cases}$$

Since $R\alpha_{v,*}\mathbb{Z}$ has trivial cohomology in degree greater than 1, we have an exact triangle $\mathbb{Z}[0] \to R\alpha_*\mathbb{Z} \to R^1\alpha_{v,*}\mathbb{Z}[-1] \to$. On the other hand, by Prop. 3.7, $H^1(U_W, \mathbb{Z}) = 0$ for all open étale subset $U$ of $\overline{Y}_{\acute{e}t}$, so $R^1\mathbb{Z} = 0$ by Theorem 1.2.

And then we have an exact triangle $\mathbb{Z}[0] \longrightarrow \tau^{\leq 3} R\mathbb{Z} \longrightarrow R^2\mathbb{Z}[-2] \longrightarrow$. Together with Lemma 4.1 and Theorem 3.2, there exists a morphism of exact triangles

$$\begin{array}{ccccccc} \mathrm{R\Gamma}(\overline{Y}_{\acute{e}t}, \mathbb{Z}) & \longrightarrow & \tau^{\leq 3}\,\mathrm{R\Gamma}(\overline{Y}_W, \mathbb{Z}) & \longrightarrow & \mathrm{R\Gamma}(\overline{Y}_{\acute{e}t}, R^2\mathbb{Z})[-2] & \longrightarrow & \\ \downarrow & & \downarrow & & \vdots\, \exists\,\gamma & & \\ \bigoplus_{v\notin Y} \mathrm{R\Gamma}(k(v)_{\acute{e}t}, \mathbb{Z}) & \longrightarrow & \bigoplus_{v\notin Y} \mathrm{R\Gamma}(k(v)_W, \mathbb{Z}) & \longrightarrow & \bigoplus_{v\notin Y} \mathrm{R\Gamma}(k(v)_{\acute{e}t}, R^1\alpha_{v,*}\mathbb{Z})[-1] & \longrightarrow & , \end{array}$$

where the rows are exact triangles, and the vertical morphisms in the first 2 columns are the canonical ones. By the remark above, we get an exact triangle

$$\mathrm{R\Gamma}_c(Y_{\acute{e}t}, \mathbb{Z}) \longrightarrow \tau^{\leq 3}\,\mathrm{R\Gamma}_c(Y_W, \mathbb{Z}) \longrightarrow \widetilde{\mathrm{R\Gamma}}_c(Y_{\acute{e}t}, R^2\mathbb{Z})[-1] \longrightarrow ,$$

where

$$\widetilde{\mathrm{R\Gamma}}_c(Y_{\acute{e}t}, R^2\mathbb{Z}) := \mathrm{cone}\Big(\mathrm{R\Gamma}(\overline{Y}_{\acute{e}t}, R^2\mathbb{Z})[-2] \xrightarrow{\gamma} \bigoplus_{v\in S} \mathrm{R\Gamma}(k(v)_{\acute{e}t}, R^1\alpha_{v,*}\mathbb{Z})[-1]\Big).$$

Since $\mathrm{R\Gamma}(\overline{Y}_{\acute{e}t}, R^2\mathbb{Z}) = \mathrm{Hom}(\mathcal{O}_F^\times, \mathbb{Q})[0]$ (Theorem 3.2) and $R^1\alpha_{v,*}\mathbb{Z} = 0$ (*resp.* $\mathbb{Q}$) when $v$ is an infinite place (resp. a finite place),

$$\widetilde{\mathrm{R\Gamma}}_c(Y_{\acute{e}t}, R^2\mathbb{Z}) = \mathrm{cone}\Big(\mathrm{Hom}(\mathcal{O}_F^\times, \mathbb{Q})[-2] \to \bigoplus_{v\in S_{<\infty}} \mathbb{Q}[-1]\Big).$$

By looking at the long exact sequence of cohomology, one can identify $\widetilde{\mathrm{R\Gamma}}_c(Y_{\acute{e}t}, R^2\mathbb{Z})$ with $\mathrm{Hom}(\mathrm{U_S}, \mathbb{Q})[-1]$. $\square$

**Theorem 4.2.** *Let* $Y = \mathrm{Spec}\,\mathcal{O}_{F,S}$*, then*

$$H^p_c(Y_W, \mathbb{Z}) = \begin{cases} 0 & p = 0, \\ \prod_S \mathbb{Z}/\mathbb{Z} & p = 1, \\ \mathrm{Pic}(Y)^{\mathcal{D}} \times \mathrm{Hom}(\mathrm{U_S}, \mathbb{Z}) & p = 2, \\ \mu_F^{\mathcal{D}} & p = 3. \end{cases}$$

*Proof.* By the definition of the cohomology with compact support, we have the following long exact sequence

$$0 \longrightarrow H^0_c(Y_W, \mathbb{Z}) \longrightarrow \mathbb{Z} \xrightarrow{\triangle} \oplus_{v\in S}\mathbb{Z} \longrightarrow H^1_c(Y_W, \mathbb{Z}) \longrightarrow 0 \to$$

$$\to \oplus_{v\in S_{<\infty}}\mathbb{Z} \longrightarrow H^2_c(Y_W, \mathbb{Z}) \longrightarrow \mathrm{Pic}^1(\overline{Y})^{\mathcal{D}} \longrightarrow 0 \longrightarrow H^3_c(Y_W, \mathbb{Z}) \longrightarrow \mu_F^{\mathcal{D}} \longrightarrow 0.$$

Since $\triangle$ is the diagonal map, we see that $H^0_c(Y_W, \mathbb{Z}) = 0$ and $H^1_c(Y_W, \mathbb{Z}) = \mathbb{Z}^{|S|}/\mathbb{Z}$. And $H^3_c(Y_W, \mathbb{Z}) = \mu_F^{\mathcal{D}}$, trivially. Finally, the exact triangle in Theorem 4.1 induces a long exact sequence of cohomology groups

$$\begin{array}{ccccccccccc} 0 \longrightarrow & H^2_c(Y_{\acute{e}t}, \mathbb{Z}) & \longrightarrow & H^2_c(Y_W, \mathbb{Z}) & \longrightarrow & \mathrm{Hom}(\mathrm{U_S}, \mathbb{Q}) & \longrightarrow & H^3_c(Y_{\acute{e}t}, \mathbb{Z}) & \longrightarrow & \mu_F^{\mathcal{D}} & \longrightarrow 0 \\ & \| & & & & & & \| & & & \\ & \mathrm{Pic}(Y)^{\mathcal{D}} & & & & & & \mathrm{Hom}(\mathrm{U_S}, \mathbb{Q}/\mathbb{Z}) & & & \end{array}$$

Since $\mathrm{Hom}(\mathrm{U_S}, \mathbb{Q})$ is torsion free, the injection on the left implies $H^2_c(Y_W, \mathbb{Z})_{tor} = \mathrm{Pic}(Y)^{\mathcal{D}}$. On the other hand, the long exact sequence in the beginning of the proof shows $\mathrm{rank}_{\mathbb{Z}}(H_c^2(Y_W,\mathbb{Z})) = \#\{v | v \notin Y\} = \mathrm{rank}_{\mathbb{Z}}(U_S)$, as $\mathrm{Pic}^1(\overline{Y})^{\mathcal{D}} = \mathrm{Pic}(\mathcal{O}_F)^{\mathcal{D}} \times \mathrm{Hom}(\mathcal{O}_F^{\times},\mathbb{Z})$([8] Prop. 6.4). As any free abelian subgroup of $\mathrm{Hom}(\mathrm{U_S},\mathbb{Q})$ of rank $\#\,\mathrm{U_S}$ is isomorphic to $\mathrm{Hom}(\mathrm{U_S},\mathbb{Z})$, Consequently,

$$H_c^2(Y_W,\mathbb{Z}) = \mathrm{Pic}(Y)^{\mathcal{D}} \times \mathrm{Hom}(\mathrm{U_S},\mathbb{Z}).$$

□

By a direct computation of the long exact sequence defining the cohomologies with compact support, we see that

**Proposition 4.1.**

$$H_c^p(Y_W,\widetilde{\mathbb{R}}) = \begin{cases} \prod_S \mathbb{R}/\mathbb{R} & p = 1,2, \\ 0 & \textit{otherwise}. \end{cases}$$

In the end of this section, we compare $\mathrm{R}\Gamma_c((\mathrm{Spec}\,\mathcal{O}_{F,S_0})_W, A)$ and $\mathrm{R}\Gamma_c((\mathrm{Spec}\,\mathcal{O}_{F,S_1})_W, A)$ for future use, where $\overline{Y} \supset S_1 \supset S_0$.

**Proposition 4.2.** *Let $\overline{Y} \supset S_1 \supset S_0 \supset \overline{Y}_\infty$ and $\mathscr{F}$ a compatible system of sheaves on the sites $(T_{K/F,S})$, then we have an exact triangle*

$$\mathrm{R}\Gamma_c((\mathrm{Spec}\,\mathcal{O}_{F,S_1})_W,\mathscr{F}) \to \mathrm{R}\Gamma_c((\mathrm{Spec}\,\mathcal{O}_{F,S_0})_W,\mathscr{F}) \to \bigoplus_{v\in S_1\setminus S_0} \mathrm{R}\Gamma(k(v)_W, i_v^*\mathscr{F}).$$

*Proof.* We first have the following commutative diagram of complexes in which the rows are exact triangles,

$$\begin{array}{ccccccc}
\mathrm{R}\Gamma_c(\mathrm{Spec}\,\mathcal{O}_{F,S_1},\mathscr{F}) & \longrightarrow & \mathrm{R}\Gamma(\overline{Y},\mathscr{F}) & \longrightarrow & \bigoplus_{v\in S_1}\mathrm{R}\Gamma(k(v)_W, i_v^*\mathscr{F}) & \longrightarrow & \\
\vdots & & \| & & \downarrow & & \\
\mathrm{R}\Gamma_c(\mathrm{Spec}\,\mathcal{O}_{F,S_0},\mathscr{F}) & \longrightarrow & \mathrm{R}\Gamma(\overline{Y},\mathscr{F}) & \longrightarrow & \bigoplus_{v\in S_0}\mathrm{R}\Gamma(k(v)_W, i_v^*\mathscr{F}) & \longrightarrow & .
\end{array}$$

This diagram can be completed by the following 3 by 3 commutative diagram of complexes,

(7)

$$\begin{array}{ccccccc}
\mathrm{R}\Gamma_c(\mathrm{Spec}\,\mathcal{O}_{F,S_1},\mathscr{F}) & \longrightarrow & \mathrm{R}\Gamma(\overline{Y},\mathscr{F}) & \longrightarrow & \bigoplus_{v\in S_1}\mathrm{R}\Gamma(k(v)_W, i_v^*\mathscr{F}) & \longrightarrow & \\
\downarrow & & \| & & \downarrow & & \\
\mathrm{R}\Gamma_c(\mathrm{Spec}\,\mathcal{O}_{F,S_0},\mathscr{F}) & \longrightarrow & \mathrm{R}\Gamma(\overline{Y},\mathscr{F}) & \longrightarrow & \bigoplus_{v\in S_0}\mathrm{R}\Gamma(k(v)_W, i_v^*\mathscr{F}) & \longrightarrow & \\
\downarrow & & \downarrow & & \downarrow & & \\
\bigoplus_{v\in S_1\setminus S_0}\mathrm{R}\Gamma(k(v)_W, i_v^*\mathscr{F}) & \longrightarrow & 0 & \longrightarrow & \bigoplus_{v\in S_1\setminus S_0}\mathrm{R}\Gamma(k(v)_W, i_v^*\mathscr{F})[1] & \longrightarrow & , \\
\downarrow & & \downarrow & & \downarrow & &
\end{array}$$

where all the rows and columns are exact triangles. Then the result follows from the first column. □

**Corollary 4.1.** *Let $A$ be $\mathbb{Z}$ or $\mathbb{R}$ on which $W_F$ acts trivially, then there is an exact triangle*

$$\tau^{\leq 3}\,\mathrm{R}\Gamma_c((\mathrm{Spec}\,\mathcal{O}_{F,S_1})_W, A) \to \tau^{\leq 3}\,\mathrm{R}\Gamma_c((\mathrm{Spec}\,\mathcal{O}_{F,S_0})_W, A) \to \bigoplus_{v\in S_1\setminus S_0}\mathrm{R}\Gamma(k(v)_W, A) \to .$$

*Proof.* This follows easily from the fact $\tau^{\geq 2}\,\mathrm{R}\Gamma(k(v)_W, A) = 0$ when $A = \mathbb{Z}, \mathbb{R}$. $\square$

**Corollary 4.2.** *For $A = \mathbb{Z}$, the long exact sequence induced by the exact triangle in the above corollary splits into short exact sequences*

$$0 \to \coprod_{S_1\setminus S_0} H^i(k(v)_W, \mathbb{Z}) \to H_c^{i+1}((\operatorname{Spec}\mathcal{O}_{F,S_1})_W, \mathbb{Z}) \to H_c^{i+1}((\operatorname{Spec}\mathcal{O}_{F,S_0})_W, \mathbb{Z}) \to 0,$$

*for all $i \geq 0$.*

*Proof.* It is suffice to prove the induced morphisms $H_c^i((\operatorname{Spec}\mathcal{O}_{F,S_0})_W, \mathbb{Z}) \to \coprod_{S_1\setminus S_0} H^i(k(v)_W, \mathbb{Z})$ are zero maps for all $i$. As $H_c^i((\operatorname{Spec}\mathcal{O}_{F,S_0})_W, \mathbb{Z})$ is concentrated in degree 1,2 and 3 and $H^i(k(v)_W, \mathbb{Z})$ is concentrated in degree 0 and 1, we only need to consider the case $i = 1$. Since the 3rd row of diagram (7) comes from the splitting exact sequence of complexes

$$0 \to \bigoplus_{v\in S_1\setminus S_0} \mathrm{R}\Gamma(k(v)_W, \mathbb{Z}) \to \bigoplus_{v\in S_1} \mathrm{R}\Gamma(k(v)_W, \mathbb{Z}) \to \bigoplus_{v\in S_0} \mathrm{R}\Gamma(k(v)_W, \mathbb{Z}) \to 0.$$

It is then clear that the morphisms

$$\bigoplus_{v\in S_0} H^{i-1}(k(v)_W, \mathbb{Z}) \xrightarrow{\delta} \bigoplus_{v\in S_1\setminus S_0} H^i(k(v)_W, \mathbb{Z})$$

are zero maps for $i \geq 0$. Therefore, the morphism of long exact sequences induced by the morphism of the 2nd and the 3rd rows in diagram (7) looks like

$$\begin{array}{ccccccccc}
\longrightarrow & \bigoplus_{v\in S_0} H^0(k(v)_W, \mathbb{Z}) & \xrightarrow{f} & H_c^1((\operatorname{Spec}\mathcal{O}_{F,S_0})_W, \mathbb{Z}) & \longrightarrow & H^1(\overline{Y}_W, \mathbb{Z})(=0) & \longrightarrow \\
 & \downarrow{\scriptstyle 0} & & \downarrow{\scriptstyle h} & & \downarrow{\scriptstyle 0} & \\
\longrightarrow & \bigoplus_{v\in S_1\setminus S_0} H^1(k(v)_W, \mathbb{Z}) & \xrightarrow{\sim} & \bigoplus_{v\in S_1\setminus S_0} H^1(k(v)_W, \mathbb{Z}) & \longrightarrow & 0 & \longrightarrow
\end{array}$$

As $f$ is surjective and the diagram is commutative, $h$ has to be the zero map. $\square$

4.2. **Verification of the Axioms of Weil-Étale Cohomology.** As before, $F$ is totally imaginary and $Y = \operatorname{Spec}\mathcal{O}_{F,S}$. In this subsection, we verify that the axioms (a)-(d) stated in the introduction hold for the generalized Weil-étale cohomology theory.

4.2.1. *Axiom (a).* We have seen that $H_c^p(Y_W, \mathbb{Z})$ are finitely generated abelian groups for $p = 0, 1, 2, 3$. Unfortunately, $H_c^p(Y_W, \mathbb{Z})$ is still of infinite rank, for odd $p \geq 5$ as $H_c^p(Y_W, \mathbb{Z}) \cong H_c^p(\operatorname{Spec}\mathcal{O}_F, \mathbb{Z})$.

4.2.2. *Axiom (b).* We need to show that $H_c^p(Y_W, \mathbb{Z}) \otimes_{\mathbb{Z}} \mathbb{R} \simeq H_c^p(Y_W, \mathbb{R})$. But, obviously, this is false for odd $p \geq 5$, as $H_c^p(Y_W, \mathbb{R}) = 0$ then. But we can prove the following.

**Proposition 4.3.** $H_c^p(Y_W, \mathbb{Z}) \otimes_{\mathbb{Z}} \mathbb{R} \simeq H_c^p(Y_W, \mathbb{R})$ *for $p = 0, 1, 2, 3$.*

*Proof.* By corollary 4.1 and taking tensor product over $\mathbb{Z}$, we get a morphism between exact triangles

$$\begin{array}{ccccccc}
\tau^{\leq 3}\,\mathrm{R}\Gamma_c((\operatorname{Spec}\mathcal{O}_{F,S})_W, \mathbb{Z}) \otimes_{\mathbb{Z}} \mathbb{R} & \longrightarrow & \tau^{\leq 3}\,\mathrm{R}\Gamma_c((\operatorname{Spec}\mathcal{O}_F)_W, \mathbb{Z}) \otimes_{\mathbb{Z}} \mathbb{R} & \longrightarrow & \bigoplus_{v\in S\setminus S_\infty} \mathrm{R}\Gamma(k(v)_W, \mathbb{Z}) \otimes_{\mathbb{Z}} \mathbb{R} & \longrightarrow \\
\downarrow & & \downarrow{\scriptstyle qis} & & \downarrow{\scriptstyle qis} & \\
\tau^{\leq 3}\,\mathrm{R}\Gamma_c((\operatorname{Spec}\mathcal{O}_{F,S})_W, \mathbb{R}) & \longrightarrow & \tau^{\leq 3}\,\mathrm{R}\Gamma_c((\operatorname{Spec}\mathcal{O}_F)_W, \mathbb{R}) & \longrightarrow & \bigoplus_{v\in S\setminus S_\infty} \mathrm{R}\Gamma(k(v)_W, \mathbb{R}) & \longrightarrow .
\end{array}$$

Clearly, for $v \nmid \infty$, $H^*(k(v)_W, \mathbb{Z}) \otimes_{\mathbb{Z}} \mathbb{R} \simeq H^*(k(v)_W, \mathbb{R})$, so the right vertical map is a quasi-isomorphism. Also, the middle one is a quasi-isomorphism ([8] Thm. 8.1). Thus, by the property of exact triangles, the left vertical morphism is also a quasi-isomorphism. □

4.2.3. *Axiom (d).* Recall that the L-function $\zeta_Y(z)$ for $Y$ is the same as the $S$-zeta function $\zeta_S(s) := \sum_{\mathfrak{a} \subset \mathcal{O}_{F,S}} N(\mathfrak{a})^{-s}$. The analytic class number formula shows its Taylor expansion at $z = 0$ is

$$\zeta_S(z) = -\frac{h_S R_S}{w} z^{|S|-1} + O(z^{|S|}),$$

where $h_S = \# \operatorname{Pic} Y$, $R_S$ the $S$-regulator, and $w = \# \mu_F$.

As $H^p_c(Y, \mathbb{Z})$ does not vanish for all even $p \geq 4$, we are unable to show that $ord_{s=0}\zeta_Y(s) = \sum_{i \geq 0} (-1)^i i \operatorname{rank}_{\mathbb{Z}} H^i_c(X, \mathbb{Z})$. However, by the computation in previous section, we see that

$$\begin{aligned} \sum_{i=0}^{3} (-1)^i i \operatorname{rank}_{\mathbb{Z}} H^i_c(X, \mathbb{Z}) &= -(|S| - 1) + 2(|S| - 1) \\ &= |S| - 1 \\ &= ord_{z=0}\zeta_Y(z). \end{aligned}$$

4.2.4. *Axiom (c) and (e).* Recall that the Euler character is defined as the following

**Definition 4.1.** *Let $A_0, \ldots, A_n$ be finitely generated abelian groups, $V_i := A_i \otimes_{\mathbb{Z}} \mathbb{R}$, and $T_i : V_i \to V_{i+1}$ are $\mathbb{R}$-linear maps s.t.*

$$0 \to V_0 \to V_1 \to \cdots \to V_n \to 0$$

*is exact. We define the* **Euler characteristic** $\chi_c(A_0, \cdots, A_n, T_0, \cdots, T_n)$ *to be*

$$\pm \prod_i |(A_i)_{tor}|^{(-1)^i} / \mathfrak{Det}(V_0, \ldots, V_n; \mathfrak{b}_0, \cdots, \mathfrak{b}_n),$$

*where $\mathfrak{b}_i = (b_{ij})_j$ is an arbitrary choice of $\mathbb{Z}$-basis of $A_i/(A_i)_{tor}$ for all $i$ and*

$$\mathfrak{Det}(V_0, \ldots, V_n; \mathfrak{b}_0, \cdots, \mathfrak{b}_n)$$

*is the image of $\otimes_i (\wedge_j b_{ij}^{(-1)^i})$ under the canonical map*

$$\otimes_i \mathfrak{Det}^{(-1)^i} V_i \to \mathbb{R},$$

*and $b_{ij}^{-1}$ is the dual of $b_{ij}$ in $\operatorname{Hom}(V_i, \mathbb{R})$.*

*Recall that for any $\mathbb{R}$-vector space $V$, $\mathfrak{Det}\, V := \wedge^{rank(V)} V$, $\mathfrak{Det}^{-1} V := \operatorname{Hom}_{\mathbb{R}}(\mathfrak{Det}\, V, \mathbb{R})$ and $\mathfrak{Det}^0 V := \mathbb{R}$.*

*Clearly, the Euler characteristic is independent to the choice of basis.*

**Example 4.1.** *Suppose $0 \to A_0 \to A_1 \to \cdots \to A_n \to 0$ is an exact sequence of finitely generated abelian groups. It induces an exact sequence of $\mathbb{R}$-vector spaces $0 \to A_0 \otimes_{\mathbb{Z}} \mathbb{R} \to A_1 \otimes_{\mathbb{Z}} \mathbb{R} \to \cdots \to A_n \otimes_{\mathbb{Z}} \mathbb{R} \to 0$. The determinant of this complex equals to $\Pi_i |(A_i)_{tor}|^{(-1)^i}$. Therefore, the Euler characteristic of this complex is $\pm 1$.*

**Example 4.2.** *Let $\{A_{i,j}\}_{i=0,..,m; j=0,..,n}$ be a set of finitely generated non-trivial abelian groups and $B_{i,j} := A_{i,j} \otimes_{\mathbb{Z}} \mathbb{R}$. Suppose that there exist exact sequences*

$$0 \to A_{i,0} \to A_{i,1} \to \cdots \to A_{i,n} \to 0, \quad \forall\, i = 0, \cdots, m.$$

*And suppose there exist $\mathbb{R}$-linear maps $T_{i,j} : B_{i,j} \to B_{i+1,j}$ for all $i = 0, \ldots, m-1$ and $j = 0, \ldots, n$, such that the following diagram commutes and all rows and columns are exact.*

$$\begin{array}{ccccccccccc}
 & & 0 & & 0 & & & & 0 & & \\
 & & \downarrow & & \downarrow & & & & \downarrow & & \\
0 & \longrightarrow & B_{0,0} & \longrightarrow & B_{0,1} & \longrightarrow & \cdots & \longrightarrow & B_{0,n} & \longrightarrow & 0 \\
 & & \downarrow{\scriptstyle T_{0,0}} & & \downarrow{\scriptstyle T_{0,1}} & & & & \downarrow{\scriptstyle T_{0,n}} & & \\
0 & \longrightarrow & B_{1,0} & \longrightarrow & B_{1,1} & \longrightarrow & \cdots & \longrightarrow & B_{1,n} & \longrightarrow & 0 \\
 & & \downarrow{\scriptstyle T_{1,0}} & & \downarrow{\scriptstyle T_{1,1}} & & & & \downarrow{\scriptstyle T_{1,n}} & & \\
 & & \vdots & & \vdots & & & & \vdots & & \\
 & & \downarrow{\scriptstyle T_{m-1,0}} & & \downarrow{\scriptstyle T_{m-1,1}} & & & & \downarrow{\scriptstyle T_{m-1,n}} & & \\
0 & \longrightarrow & B_{m,0} & \longrightarrow & B_{m,1} & \longrightarrow & \cdots & \longrightarrow & B_{m,n} & \longrightarrow & 0 \\
 & & \downarrow & & \downarrow & & & & \downarrow & & \\
 & & 0 & & 0 & & & & 0 & &
\end{array}$$

*Then $\Pi_i\, \chi_c\big((A_{k,i})_k, (T_{k,i})_k\big)^{(-1)^i} = \pm 1$.*

*This can be proved by showing the diagram*

$$\begin{array}{ccc}
\otimes_i(\otimes_j \mathfrak{Det}^{(-1)^{i+j}} B_{ij}) & \xrightarrow{\sim} & \otimes_j(\otimes_i \mathfrak{Det}^{(-1)^{i+j}} B_{ij}) \\
 & {\scriptstyle\sim}\searrow \quad \swarrow{\scriptstyle\sim} & \\
 & \mathbb{R} &
\end{array}$$

*is commutative and all arrows are the canonical isomorphisms.(cf.*[7]*).*

Let $\mathscr{F} = ((F_{L,S}); (f_t))$ be a compatible system of abelian sheaves on the sites $(T_{L/F,S})$. And observe that $H^p(T_{L/F,S}, F_{L/F,S}) = Ext^p_{\widetilde{T_{L/F,S}}}(\mathbb{R}, F_{L/F,S})$ for any sheaf of $\widetilde{\mathbb{R}}$-modules $F$, so $H^p(\overline{Y}_W, \mathscr{F}) = \varinjlim Ext^p_{\widetilde{T_{L/F,S}}}(\mathbb{R}, F_{L,S})$. This allow us to define a cup product

$$H^p(\overline{Y}_W, \mathscr{F}) \times H^q(\overline{Y}_W, \mathbb{R}) \to H^{p+q}(\overline{Y}_W, \mathscr{F}).$$

Consider the compatible system $((i_{v,*}i_v^* F_{L,S}); (i_{v,*}i_v^* f_t))$, and take mapping cone, we can actually define a cup product

$$H_c^p(Y_W, \mathscr{F}) \times H^q(\overline{Y}_W, \mathbb{R}) \to H_c^{p+q}(Y_W, \mathscr{F}).$$

The Leray spectral sequence for $j_{K/F,S}$ yields

$$0 \to H^1(\overline{Y}_W, \mathbb{R}) \to H^1(W_F, \mathbb{R}) \to H^0(\overline{Y}_W, R^1 j_* \mathbb{R}) = 0 \text{ (see [8])}.$$

Since $H^1(W_F, \mathbb{R}) = \mathrm{Hom}_{cont}(W_F, \mathbb{R}) \cong \mathrm{Hom}_{cont}(\mathbb{R}, \mathbb{R})$([8]Lemma 3.4) , we may choose an element $\psi$ of $H^1(\overline{Y}_W, \mathbb{R})$ which corresponds to the composition

$$W_F \xrightarrow{ab} C_F \xrightarrow{\log|.|} \mathbb{R}.$$

Note that $\psi$ corresponds to the identity map in $\mathrm{Hom}_{cont}(\mathbb{R},\mathbb{R})$, and $\psi \smile \psi = 0$, so the sequence $(H^i_c(Y,\mathscr{F}), \smile \psi)$:

$$\cdots \to H^i_c(Y_W,\mathscr{F}) \xrightarrow{\smile\psi} H^{i+1}_c(Y_W,\mathscr{F}) \to \cdots$$

is a complex.

In the following, we will prove that $(H^i_c(Y_W,\mathbb{R}), \smile \psi)$ is exact. Together with the isomorphisms $H^i_c(Y_W,\widetilde{\mathbb{R}}) \cong H^i_c(Y_W,\mathbb{Z}) \otimes_{\mathbb{Z}} \mathbb{R}$ for $i \leq 3$, the Euler characteristic $\chi_c\big((H^q_c(Y_W,\mathbb{Z}))_{q\leq 3}, (\smile \psi)\big)$ is well-defined and called the Euler characteristic $\chi_c(Y)$ for $Y$.

It was already known that $\chi_c(\operatorname{Spec}\mathcal{O}_F) = Rh/w$.([8] Thm. 8.1).

**Proposition 4.4.** (a) *The complex $(H^i_c(Y_W,\widetilde{\mathbb{R}}), \smile \psi)$ is exact.*

(b) *Moreover, its Euler character equals $\pm\frac{R_S h_S}{w}$, where $R_S$ is the $S$-regulator.*

*Proof.* (a) Recall that $H^i_c(Y_W,\widetilde{\mathbb{R}}) \neq 0$ only for $i = 1, 2$. So we only need to show

$$H^1_c(Y_W,\widetilde{\mathbb{R}}) \xrightarrow{\smile\psi} H^2_c(Y_W,\widetilde{\mathbb{R}})$$

is an isomorphism. For this, consider the morphism between two exact sequences

(8)
$$\begin{array}{ccccccccc} 0 & \to & \coprod_{S\setminus S_\infty} H^0(k(v)_W,\mathbb{R}) & \to & H^1_c(Y_W,\widetilde{\mathbb{R}}) & \to & H^1_c((\operatorname{Spec}\mathcal{O}_F)_W,\widetilde{\mathbb{R}}) & \to & 0 \\ & & \downarrow{\smile\psi} & & \downarrow{\smile\psi} & & \downarrow{\smile\psi} & & \\ 0 & \to & \coprod_{S\setminus S_\infty} H^1(k(v)_W,\mathbb{R}) & \to & H^2_c(Y_W,\widetilde{\mathbb{R}}) & \to & H^2_c((\operatorname{Spec}\mathcal{O}_F)_W,\widetilde{\mathbb{R}}) & \to & 0, \end{array}$$

in which the rows are obtained by tensoring with $\mathbb{R}$ over $\mathbb{Z}$ the exact sequences derived in corollary 4.2. Recall that the right vertical map is an isomorphism ([8]Thm. 8.1). We are now showing the left vertical map is also an isomorphism. Note that $i_v$ is an embedding and preserves injectives, so that we have the following commutative diagram

$$\begin{array}{ccccc} H^0(k(v)_W,\mathbb{R}) & \times & H^1(\overline{Y},\widetilde{\mathbb{R}}) & \xrightarrow{\smile} & H^1(k(v)_W,\mathbb{R}) \\ \| & & \downarrow & & \| \\ H^0(k(v)_W,\mathbb{R}) & \times & H^1(k(v)_W,\mathbb{R}) & \xrightarrow{\smile} & H^1(k(v)_W,\mathbb{R}). \end{array}$$

For any finite place in $S$, the diagram (2) ensures that the image of $\psi$ in $H^1(k(v)_W,\mathbb{R}) = Hom(W_{k(v)},\mathbb{R})$ is $l_v : \sigma \mapsto \log N(v)$, where $\sigma$ is the generator of $W_{k(v)}$. Let $1_v$ is the identity element in $H^0(k(v)_W,\mathbb{R}) = \mathbb{R}$. Then $1_v \smile \psi = 1_v \smile l_v = l_v$. So, cup product with $\psi$ identifies $\coprod_{S\setminus S_\infty} H^0(k(v)_W,\mathbb{R})$ and $\coprod_{S\setminus S_\infty} H^1(k(v)_W,\mathbb{R})$. Therefore, the middle vertical map of diagram (8) is also an isomorphism.

(b) The rows of (8) are induced by tensoring with $\mathbb{R}$ the exact sequences

$$0 \to \coprod_{S_{<\infty}} H^i(k(v)_W,\mathbb{Z}) \to H^{i+1}_c(Y_W,\mathbb{Z}) \to H^{i+1}_c((\operatorname{Spec}\mathcal{O}_F)_W,\mathbb{Z}) \to 0, \;\; i = 0, 1.$$

Let $\chi_0 := \chi_c\Big((0, H^0(k(v)_W,\mathbb{Z}), H^1(k(v)_W,\mathbb{Z})), (\smile \psi)\Big)$.

Applying Example 4.2 to (8), we see that

$$\chi_c(Y) = \pm\chi_c(\operatorname{Spec}\mathcal{O}_F)\chi_0.$$

Because $H^i(k(v)_W, \mathbb{Z})$ is torsion free, and we has seen in (a) that $\smile\psi$ sends $1_v$ to $l_v : \sigma \mapsto \log N(v)$, so $\chi_0 = \pm\Pi_{S_{<\infty}} \log N(v)$. Thus,

$$\chi_c(Y) = \pm\frac{Rh}{w}\prod_{S_{<\infty}} \log N(v) = \pm\frac{R_S h_S}{w}$$

(c.f. [14] Lemma 2.1 for the last equality).

□

4.3. **A Canonical Representation of Tate Sequences.** As before $Y = \operatorname{Spec}\mathcal{O}_{F,S}$ and $G = \operatorname{Gal}(F/L)$ for a subfield $L$ of $F$. Recall that

$$H^p_c(Y_W, \mathbb{Z}) = \begin{cases} 0 & p = 0, \\ \prod_S \mathbb{Z}/\mathbb{Z}\ (= X_S^\vee) & p = 1, \\ \operatorname{Hom}(\mathrm{U_S}, \mathbb{Z}) \oplus \operatorname{Pic}(Y)^D & p = 2, \\ \mu_F^D & p = 3. \end{cases}$$

**Proposition 4.5.** *The cohomologies of* $\operatorname{RHom}_{\mathbb{Z}}(\tau^{\leq 3}\operatorname{R\Gamma}_c(Y_W, \mathbb{Z}), \mathbb{Z})[-2]$ *concentrate at degrees 0 and 1 and* $H^0 = U_S$*,* $H^1 = X_S \oplus \operatorname{Pic}(Y)$*.*

*Proof.* Consider the spectral sequence

$$E_2^{pq} = \operatorname{Ext}^p_{\mathbb{Z}}(H^{-q}(C^\bullet), \mathbb{Z}) \Rightarrow H^{p+q}(\operatorname{RHom}(C^\bullet, \mathbb{Z})),$$

in which $C^\bullet = \tau^{\leq 3}\operatorname{R\Gamma}_c(Y_W, \mathbb{Z})$.

It is clear that $E_2^{pq} = 0$ for all $p < 0$ or $q > 0$, i.e. the $E_2$ has only non-trivial terms on the 4th quadrant

$$\begin{array}{ccc} 0 & 0 & \\ X_S = X_S^{\vee\vee} & 0 & \\ U_S/tor = U_S^{\vee\vee} & \operatorname{Ext}^1((\operatorname{Pic}Y)^D, \mathbb{Z}) & \\ 0 & \operatorname{Ext}^1(\mu_F^D, \mathbb{Z}) & 0 \end{array}$$

We have exact sequences

$$0 \to \operatorname{Pic}Y \to H^{-1} \to X_S \to 0 \quad \text{and} \quad 0 \to \mu_F \to H^{-2} \to U_S/tor \to 0.$$

Thus, $H^{-1} = X_S \oplus \operatorname{Pic}(Y)$ and $H^{-2} = (U_S/tor) \oplus \mu_F = U_S$. This shows that $\operatorname{RHom}(\tau^{\leq 3}\operatorname{R\Gamma}_c(Y_W, \mathbb{Z}), \mathbb{Z})[-2]$ has only two non-trivial cohomologies $H^0 = U_S$ and $H^1 = X_S \oplus \operatorname{Pic}Y$. □

From now on we assume that $S$ is large enough so that $\operatorname{Pic}Y = 0$ and is stable under the action of $G$. We will prove $\operatorname{RHom}(\tau^{\leq 3}\operatorname{R\Gamma}_c(Y_W, \mathbb{Z}), \mathbb{Z})[-2]$ is indeed quasi-isomorphic to $\Psi_S$, the canonical Tate sequence. (see the introduction section).

**Lemma 4.2.** *Let* $\mathcal{A}$ *and* $\mathcal{B} \in \mathcal{D}$*, the derived category of* $\mathbb{Z}[G]$*-modules.*

a) *If* $H^i(\mathcal{B})$ *are* $\mathbb{Z}[G]$*-injective for all* $i$*, then*

$$\operatorname{Hom}_{\mathcal{D}}(\mathcal{A}, \mathcal{B}) \cong \prod_i \operatorname{Hom}_G(H^i(\mathcal{A}), H^i(\mathcal{B})).$$

b) *If* $H^i(\mathcal{A})$ *are* $\mathbb{Z}[G]$*-cohomologically trivial, then we have an exact sequence*

$$0 \to \prod_i \operatorname{Ext}^1_G(H^i(\mathcal{A}), H^{i-1}(\mathcal{B})) \to \operatorname{Hom}_{\mathcal{D}}(\mathcal{A}, \mathcal{B}) \to \prod_i \operatorname{Hom}_G(H^i(\mathcal{A}), H^i(\mathcal{B})) \to 0.$$

*Proof.* a) By [15], there is a spectral sequence

$$\prod_i \mathrm{Ext}^p_G(H^i(\mathcal{A}), H^{q+i}(\mathcal{B})) \Rightarrow H^{p+q}(\mathrm{RHom}_{\mathbb{Z}[G]}(\mathcal{A},\mathcal{B})).$$

By condition, $\mathrm{Ext}^p(H^i(\mathcal{A}), H^{q+i}(\mathcal{B}))$ is non-trivial only if $p = 0$, so the above spectral sequence simply shows

$$\mathrm{Hom}_{\mathcal{D}}(\mathcal{A},\mathcal{B}) = H^0(\mathrm{RHom}_{\mathbb{Z}[G]}(\mathcal{A},\mathcal{B})) = \prod_i \mathrm{Hom}_G(H^i(\mathcal{A}), H^i(\mathcal{B})).$$

b) In the same spectral sequence, since $H^i(\mathcal{A})$ are cohomologically trivial, $H^i(\mathcal{A})$ are of projective dimension 1, so $\mathrm{Ext}^p_G(H^i(\mathcal{A}), H^{q+i}(\mathcal{B})) = 0$ for all $p \geq 2$. Therefore, on the position with $p + q = 0$, we have only 2 non-trivial groups $E_2^{0,0}$ and $E_2^{1,-1}$, and both the difference maps are 0. This gives the desired short exact sequence.

□

**Lemma 4.3.** *The canonical exact triangle*

$$\mathrm{RHom}_{\mathbb{Z}}(\mathrm{R}\Gamma_c(Y_{\acute{e}t},\mathbb{Z}),\mathbb{Z}) \to \mathrm{RHom}_{\mathbb{Z}}(\mathrm{R}\Gamma_c(Y_{\acute{e}t},\mathbb{Z}),\mathbb{Q}) \to \mathrm{RHom}_{\mathbb{Z}}(\mathrm{R}\Gamma_c(Y_{\acute{e}t},\mathbb{Z}),\mathbb{Q}/\mathbb{Z}) \to$$

*is isomorphic to the exact triangle*

$$\Psi'_S[2] \to X_S \otimes_{\mathbb{Z}} \mathbb{Q}[1] \to \widetilde{\Psi_S}'[3] \to$$

*where $\Psi'_S$ is the image of $\Psi_S$ in $\mathrm{Ext}^2_G(X_S, \widehat{U}_S)$ and $\widetilde{\Psi_S}'$ is the image of $\Psi_S$ in $\mathrm{Ext}^3_G(X_S \otimes_{\mathbb{Z}} \mathbb{Q}/\mathbb{Z}, \widehat{U}_S)$.*

*Proof.* Suppose $\Psi_S$ is the complex $A \to B$ and $\widetilde{\Psi}_S$ is the complex $A \to B \to X_S \otimes_{\mathbb{Z}} \mathbb{Q}$. Recall that $\Psi_S$ presents a class of $\mathrm{Ext}^2_G(X_S, \mathrm{U_S})$ and $\Psi'_S$ the image of $\Psi_S$ under the map $\mathrm{Ext}^2_G(X_S, \mathrm{U_S}) \xrightarrow{\sim} \mathrm{Ext}^2_G(X_S, \widehat{U}_S)$. Indeed , $\Psi'_S$ is the complex $(A \oplus \hat{U}_S)/U_S =: A' \to B$. Similarly, $\widetilde{\Psi}'_S$ is the complex $A' \to B \to X_S \otimes_{\mathbb{Z}} \mathbb{Q}$. Recall that $\mathrm{RHom}(\mathrm{R}\Gamma_c(Y_{\acute{e}t},\mathbb{Z}),\mathbb{Q}/\mathbb{Z}) \xrightarrow{\text{qis}} \widetilde{\Psi}'_S$ by definition ([3] Prop. 3.1). Also, as $\mathbb{Q}$ is injective, $H^i(\mathrm{RHom}_{\mathbb{Z}}(\mathrm{R}\Gamma_c(Y_{\acute{e}t},\mathbb{Z}),\mathbb{Q})) = \mathrm{Hom}_{\mathbb{Z}}(H_c^{-i}(Y_{\acute{e}t},\mathbb{Z}),\mathbb{Q})$. Thus,

$$\mathrm{RHom}_{\mathbb{Z}}(\mathrm{R}\Gamma_c(Y_{\acute{e}t},\mathbb{Z}),\mathbb{Q}) \xrightarrow{\text{qis}} X_S \otimes_{\mathbb{Z}} \mathbb{Q}[1].$$

We claim that the following diagram commutes for a suitable choice of the left vertical quasi-isomorphism.

(9)
$$\begin{array}{ccc} \mathrm{RHom}_{\mathbb{Z}}(\mathrm{R}\Gamma_c(Y_{\acute{e}t},\mathbb{Z}),\mathbb{Q})[-3] & \longrightarrow & \mathrm{RHom}_{\mathbb{Z}}(\mathrm{R}\Gamma_c(Y_{\acute{e}t},\mathbb{Z}),\mathbb{Q}/\mathbb{Z})[-3] \\ {\scriptstyle qis}\uparrow & & \uparrow{\scriptstyle qis} \\ X_S \otimes_{\mathbb{Z}} \mathbb{Q}[-2] & \longrightarrow & \widetilde{\Psi}'_S. \end{array}$$

Since the both horizontal maps induce the canonical projection $X_S \otimes_{\mathbb{Z}} \mathbb{Q} \to X_S \otimes_{\mathbb{Z}} \mathbb{Q}/\mathbb{Z}$ on $H^2$, the right vertical quasi-isomorphism induces the identity map on $H^2$ and the left vertical one can be induced by any $G$-automorphism of $X_S \otimes_{\mathbb{Z}} \mathbb{Q}$, so there is a suitable choice of the quasi-isomorphism so that the above diagram commutes for $H^2$ groups.

Let $\mathcal{A}$ be $X_S \otimes_{\mathbb{Z}} \mathbb{Q}[-2]$ and $\mathcal{B}$ be $\mathrm{RHom}_{\mathbb{Z}}(\mathrm{R}\Gamma_c(Y_{\acute{e}t},\mathbb{Z}),\mathbb{Q}/\mathbb{Z})[-3]$. Clearly, $H^i(\mathcal{A})$ are $\mathbb{Z}[G]$-cohomologically trivial. So we can apply lemma 4.2(b) to $\mathcal{A}$ and $\mathcal{B}$, and it induces an exact sequence

$$0 = \mathrm{Ext}^1_G(X_S \otimes_{\mathbb{Z}} \mathbb{Q}, H^1(\mathcal{B})(=0)) \longrightarrow \mathrm{Hom}_{\mathcal{D}}(\mathcal{A},\mathcal{B}) \longrightarrow \mathrm{Hom}_G(H^2(\mathcal{A}), H^2(\mathcal{B})) \to 0.$$

Thus, $\mathrm{Hom}_{\mathcal{D}}(\mathcal{A},\mathcal{B}) \cong \mathrm{Hom}_G(H^2(\mathcal{A}), H^2(\mathcal{B}))$. Therefore, the commutative on $H^2$ implies the diagram (9) indeed commutes. By taking exact triangles of each row, and note that $\mathrm{cone}(X_S \otimes_{\mathbb{Z}} \mathbb{Q}[-2] \to \widetilde{\Psi}'_S)[-1] \simeq \Psi'_S[2]$, we get the desired isomorphism of exact triangles. □

**Theorem 4.3.** *Let $Y = \operatorname{Spec}\mathcal{O}_{F,S}$, where $S$ is large enough so that $\operatorname{Pic}Y = 0$ and is stable under the action of $G$. Then*

$$\mathrm{RHom}_{\mathbb{Z}}(\tau^{\leq 3}\,\mathrm{R}\Gamma_c(Y_W,\mathbb{Z}),\mathbb{Z})[-2] \xrightarrow{\text{qis}} \Psi_S,$$

*where $\Psi_S$ is the Tate sequence representing the canonical class of $\mathrm{Ext}^2_G(X_S, U_S)$.*

*Proof.* Here, we use the same notation with previous lemma. Clearly, we have a canonical exact triangle

$$\Psi_S \to \Psi'_S \to \widehat{U}_S/U_S[0] \to .$$

On the other hand, by applying $\mathrm{RHom}_{\mathbb{Z}}(-,\mathbb{Z})$ to exact triangle (6), we have another exact triangle

$$\mathrm{RHom}_{\mathbb{Z}}(\tau^{\leq 3}\,\mathrm{R}\Gamma_c(Y_W,\mathbb{Z}),\mathbb{Z})[-2] \to \mathrm{RHom}_{\mathbb{Z}}(\mathrm{R}\Gamma_c(Y_{\acute{e}t},\mathbb{Z}),\mathbb{Z})[-2] \to \mathrm{RHom}_{\mathbb{Z}}(\mathrm{Hom}_{\mathbb{Z}}(U_S,\mathbb{Q}),\mathbb{Z})[-1]$$

The previous lemma ensures that $\Psi'_S \xrightarrow{\text{qis}} \mathrm{RHom}_{\mathbb{Z}}(\mathrm{R}\Gamma_c(Y_{\acute{e}t},\mathbb{Z}),\mathbb{Z})[-2]$.

Also, observe that $\mathrm{RHom}_{\mathbb{Z}}(\mathrm{Hom}_{\mathbb{Z}}(U_S,\mathbb{Q}),\mathbb{Z})[-1] \xrightarrow{\text{qis}} \widehat{U}_S/U_S[0]$ and we have seen that $\mathrm{RHom}_{\mathbb{Z}}(\tau^{\leq 3}\,\mathrm{R}\Gamma_c(Y_W,\mathbb{Z}),\mathbb{Z})[-2]$ has the same cohomologies as $\Psi_S$. We claim that there is an isomorphism between these two exact triangles,

$$\begin{array}{ccccccc}
\mathrm{RHom}_{\mathbb{Z}}(\tau^{\leq 3}\,\mathrm{R}\Gamma_c(Y_W,\mathbb{Z}),\mathbb{Z})[-2] & \longrightarrow & \Psi'_S & \xrightarrow{\beta} & \mathrm{RHom}_{\mathbb{Z}}(\mathrm{Hom}_{\mathbb{Z}}(U_S,\mathbb{Q}),\mathbb{Z})[-1] & \longrightarrow & \\
\uparrow{\scriptstyle\delta} & & \| & & \uparrow{\scriptstyle\alpha} & & \\
\Psi_S & \longrightarrow & \Psi'_S & \xrightarrow{\gamma} & \widehat{U}_S/U_S[0] & \longrightarrow & .
\end{array}$$

In the top exact triangle, we replaced $\mathrm{RHom}_{\mathbb{Z}}(\mathrm{R}\Gamma_c(Y_{\acute{e}t},\mathbb{Z}),\mathbb{Z})[-2]$ by $\Psi'_S$. By taking the long exact sequence induced by the top exact triangle, we see that $H^2(\beta) = i \circ p$ where $p : \widehat{U}_S \to \widehat{U}_S/U_S$ is the canonical projection and $i \in Aut(\widehat{U}_S/U_S)$. Clearly, $H^2(\gamma) = p$. As $\mathrm{RHom}_{\mathbb{Z}}(\mathrm{Hom}_{\mathbb{Z}}(U_S,\mathbb{Q}),\mathbb{Z})[-1] \xrightarrow{\text{qis}} \widehat{U}_S/U_S[0]$, we may take $\alpha$ to be the morphism induced by $i$. Then the right square commutes on $H^2$. By lemma 4.2, because $\widehat{U}_S/U_S$ is uniquely divisible and thus injective, $H^2(\beta) = H^2(\alpha) \circ H^2(\gamma)$ implies $\beta = \alpha \circ \gamma$. As $\alpha$ is a quasi-isomorphism, so is $\delta$ by property of derived categories. □

**Remark.** *It is equivalent to say that $\tau^{\leq 3}\,\mathrm{R}\Gamma_c(Y_W,\mathbb{Z}) \xrightarrow{\text{qis}} \mathrm{RHom}(\Psi_S,\mathbb{Z})[-2]$ and this is an evidence that, at least for open subschemes of spectra of number rings, the Weil-étale cohomology can be obtained from usual cohomology theories without using the Weil groups.*

## 5. The Construction of $\mathbf{R}\mathbb{G}_m$

Throughout this section, we assume that $F$ is a *totally imaginary* number field, $X = \operatorname{Spec}\mathcal{O}_F$, $U = \operatorname{Spec}\mathcal{O}_{L,S}$ be any connected etale neighborhood of $X$. We also assume that $K$ is a subfield of $L$ such that the extension $L/K$ is Galois and $G := \mathrm{Gal}(L/K)$.

When $Pic(U) = 0$ and $S$ is a $G$-stable set containing all those ramified primes, there exists a canonical Tate sequence:

$$\Psi_U : \Psi_S^0 \to \Psi_S^1,$$

which represents a canonical class of $\mathrm{Ext}_G^2(X_S, U_S)$. By the construction of Tate sequences, if $U' = \operatorname{Spec} \mathcal{O}_{L',S'}$ is étale over $U$, then there is a morphism from $\Psi_S$ to $\Psi_{S'}$ that fits into the following morphism of complexes

$$\begin{array}{ccccccccc} 0 & \longrightarrow & U_S & \longrightarrow & \Psi_S^0 & \longrightarrow & \Psi_S^1 & \longrightarrow & X_S & \longrightarrow & 0 \\ & & \downarrow & & \downarrow & & \downarrow & & \downarrow{\scriptstyle \beta} & & \\ 0 & \longrightarrow & U_{S'} & \longrightarrow & \Psi_{S'}^0 & \longrightarrow & \Psi_{S'}^1 & \longrightarrow & X_{S'} & \longrightarrow & 0, \end{array}$$

where $\beta((a_v)_{v\in S}) = \big((a_w := [L_w : F_v] \cdot a_v)_{w|v, w\in S'}\big)$.

In fact $\{X_S\}$ and $\{\beta\}$ define an $\overline{X}_{\acute{e}t}$-presheaf $\mathscr{X}$, and we denote its associated sheaf by $\widetilde{\mathscr{X}}$. In section 5.1, we construct a complex $R\,\mathbb{G}_\mathrm{m} \in Sh(\overline{X}_{\acute{e}t})$ which represents the canonical class of $\mathrm{Ext}^2_{\overline{X}_{et}}(\widetilde{\mathscr{X}}, \mathbb{G}_\mathrm{m}) = \widehat{\mathbb{Z}}$ (see below). The cohomologies of $\mathrm{R}\Gamma(U_{\acute{e}t}, R\,\mathbb{G}_\mathrm{m})$ and $\mathrm{RHom}_{\mathbb{Z}}\left(\tau^{\leq 3}\,\mathrm{R}\Gamma_c(U_W, \mathbb{Z}),\, \mathbb{Z}\right)[-2]$ are the same for *arbitrary* $U$. Furthermore, in section 4.2, we prove that

$$\mathrm{R}\Gamma(U_{\acute{e}t}, R\,\mathbb{G}_\mathrm{m}) \xrightarrow{\mathrm{qis}} \mathrm{RHom}\left(\tau^{\leq 3}\,\mathrm{R}\Gamma_c(U_W, \mathbb{Z}), \mathbb{Z}\right)[-2], \tag{10}$$

in $D(\mathbb{Z}[G])$, when $S$ is stable under the action of $G$.

This implies $\mathrm{R}\Gamma(U_{\acute{e}t}, R\,\mathbb{G}_\mathrm{m})$ and $\mathrm{RHom}_{\mathbb{Z}}\left(\tau^{\leq 3}\,\mathrm{R}\Gamma_c(U_W, \mathbb{Z}),\, \mathbb{Z}\right)[-2]$ are canonically quasi-isomorphic in $D(\mathbb{Z})$ for any $U$. Also, when $U$ is small enough so that the Tate sequence exists, then by (10) and by the help of Theorem 4.3, one can see that

$$\mathrm{R}\Gamma(U_{\acute{e}t}, R\,\mathbb{G}_\mathrm{m}) \simeq \Psi_U$$

in $D(\mathbb{Z}[G])$.

Thus, one can conclude that the complex $\mathrm{R}\Gamma(U_{\acute{e}t}, R\,\mathbb{G}_\mathrm{m})$ generalizes the Tate sequences and its $\mathbb{Z}$-dual defines Weil-etale cohomology of $S$-integers without truncations.

5.1. **The Definition.** Let $\mathscr{M}$ be the sheaf sending a connected étale $U \to \overline{X}$ to the group $\bigoplus_{w\in\overline{U}} A_w$, where $A_w$ is $\mathbb{Q}$ (resp. $\mathbb{Z}$) if $w \nmid \infty$ (resp. $w \mid \infty$), and the transition map $\mathscr{M}(U) \to \mathscr{M}(V)$ for an $\overline{X}$-morphism $V \to U$ is

$$\begin{array}{ccc} \displaystyle\bigoplus_{v\in\overline{U}} A_v & \longrightarrow & \displaystyle\bigoplus_{w\in\overline{V}} A_w \\ a_v & \mapsto & \displaystyle\sum_{w|v} [K(V)_w : K(U)_v] a_v. \end{array}$$

It is clear that the projections $\bigoplus_{v\in\overline{U}} A_v \twoheadrightarrow \bigoplus_{v\in U} A_v$, for any étale $U \to \overline{X}$, define an epimorphism of sheaves $\mathscr{M} \to \bigoplus_{v\in\overline{X}} i_{v,*} A_v$.

We define $\mathscr{F}^\bullet$ to be the complex

$$\begin{array}{ccc} \mathscr{F}^0 & \xrightarrow{\Sigma} & \mathscr{F}^1 \\ & & \| \\ & & \mathbb{Q}, \end{array}$$

where $\mathscr{F}^0 = ker(\mathscr{M} \to \bigoplus_{v \in \overline{X}} i_{v,*}A_v)$ is the $\overline{X}$-étale sheaf

$$U \mapsto \bigoplus_{w \in \overline{U} - U} A_w,$$

for any connected $U \to \overline{X}$ and $\sum(U)((a_w)) = \sum_w \frac{1}{[K(U)_w : F_v]} a_w$.

Note that if $U$ doesn't contain all the finite places of $F$, then $\sum(U)$ is surjective. This implies $\sum$ is an epimorphism in the category $Sh(\overline{X}_{ét})$ as any $U \to \overline{X}$ can be covered by some $\{U_i\}$, where $U_i$ doesn't contain all the finite places of $K(U_i)$, for all $i$. Consequently, $\mathscr{F}^\bullet \simeq (\ker \sum)[0]$ in $D(Sh(\overline{X}_{ét}))$.

**Remark.** (a) *The way that we define the transition maps of the sheaf $\mathscr{F}^0$ makes $\mathscr{X}$ a sub-presheaf of $\ker \sum$. In fact, we have the following exact sequence of presheaves*

$$0 \to \mathscr{X} \to \ker \sum \to \widetilde{Br} \to 0,$$

*where $\widetilde{Br}$ is the presheaf $U \mapsto H^2(U, \mathbb{G}_{\mathrm{m}})$. By taking associated sheaves, we see that $\widetilde{\mathscr{X}} \cong \ker \sum$ as ${}^a(H^i(-, \mathscr{G})) = 0$ for any sheaf $\mathscr{G}$ when $i \geq 1$.*

(b) *Let $M_v = \varinjlim_L M_{L,v}$ where $L$ runs over all finite extension of $F$ and $M_{L,v} = \bigoplus_w A_w$ where the sum runs over all places $w$ of $L$ lying over $v$ and the transition map is defined similarly to those of $\mathscr{X}$. We denote $\mathscr{M}_v$ the sheaf associated to the Galois module $M_v$. Equivalently, in the level of sheaves, $\mathscr{M}_v = \varinjlim_L \pi_{L,*}A_v$ where $\pi_L : \operatorname{Spec} L^v \longrightarrow \operatorname{Spec} F$ is the canonical morphism and $L^v$ is the fixed field of $D_v \cap \operatorname{Gal}(L/F)$. Clearly, $\mathscr{M} = \bigoplus_{v \in \overline{X}} \alpha_* \mathscr{M}_v$, where $\alpha : \operatorname{Spec} F \to \overline{X}$. Moreover, consider the Cartesian diagram of schemes*

$$\begin{array}{ccc} \operatorname{Spec} L^v & \xrightarrow{\alpha_L} & u(L^v) \\ \Big\downarrow{\scriptstyle \pi_L} & & \Big\downarrow{\scriptstyle \pi'_L} \\ \operatorname{Spec} F & \xrightarrow{\alpha} & \overline{X}, \end{array}$$

*where $u(L^v) = \overline{\operatorname{Spec} \mathcal{O}_{L^v}}$. Hence $\alpha_* \varinjlim_L \pi_{L,*}A_v = \varinjlim_L \alpha_* \pi_{L,*}A_v = \varinjlim_L \pi'_{L,*}\alpha_{L,*}A_v = \varinjlim_L \pi'_{L,*}A_v$. (Note that $\alpha_{L,*}A_v$ is the constant $A_v$ on $\operatorname{Spec}(L^v)$.)*

**Lemma 5.1.** *The sheaf $\mathscr{F}^0$ is cohomologically trivial.*

*Proof.* Since $\mathscr{F}^0$ fits into the short exact sequence

$$0 \to \mathscr{F}^0 \to \bigoplus_{v \in \overline{X}} \alpha_* \mathscr{M}_v \to \bigoplus_{v \in \overline{X}} i_{v,*}A_v \to 0, \tag{11}$$

which is derived from the short exact sequence of sections:

$$0 \to \bigoplus_{v \in \overline{U} - U} A_v \to \bigoplus_{v \in \overline{U}} A_v \to \bigoplus_{v \in U} A_v \to 0.$$

Since there is no real place and $A_v = \mathbb{Q}$ when $v \nmid \infty$, the cohomology of $i_{v,*}A_v$ concentrates at degree zero. Also, the above exact sequence remains exact when we pass to global sections. Thus, we only need to show that $\alpha_* \mathscr{M}_v$ is cohomologically trivial for any $v$.

Recall that $\alpha_*\mathscr{M}_v = \varinjlim_L \pi'_{L,*}A_v$. When $v$ is finite, $\alpha_* M_v$ is acyclic as $A_v = \mathbb{Q}$ and inductive limit and $\pi'_{L,*}$ both preserve acyclic sheaves (as $\pi'_L$ is finite).

When $v|\infty$, $I_v$ is trivial, so $M_v = Ind_1^{G_F} A_v$ and then $M_v$ is cohomologically trivial. Therefore $i_w^*(R^q\alpha_* M_v) = H^q(I_w, M_v) = 0$ for $q > 0$. The Leray spectral sequence for $\alpha$,

$$H^p(U_{\acute{e}t}, R^q\alpha_* M_v) \Rightarrow H^{p+q}(F, M_v),$$

for any étale $U$ over $\overline{X}$, is degenerated, i.e. $H^p(U_{\acute{e}t}, \alpha_* M_v) = H^p(F, M_v) = 0$ for $p > 0$. Thus the result follows. □

**Proposition 5.1.** *For any connected étale $U \to \overline{X}$,*

$$H^p(U_{\acute{e}t}, \mathscr{F}^\bullet) = \begin{cases} \ker \Sigma(U) & p = 0, \\ \operatorname{coker} \Sigma(U) = \mathbb{Q}/Im(\Sigma(U)) & p = 1, \\ 0 & p \geq 2. \end{cases}$$

*In particular, $H^0(\overline{X}_{\acute{e}t}, \mathscr{F}^\bullet) = 0$, $H^1(\overline{X}_{\acute{e}t}, \mathscr{F}^\bullet) = \mathbb{Q}$, $H^0((\overline{X} - v)_{\acute{e}t}, \mathscr{F}^\bullet) = 0$ and $H^1((\overline{X} - v)_{\acute{e}t}, \mathscr{F}^\bullet) = \mathbb{Q}/A_v$.*

**Corollary 5.1.** *The canonical morphism $H_v^1(\overline{X}_{\acute{e}t}, \mathscr{F}^\bullet) \to H^1(\overline{X}_{\acute{e}t}, \mathscr{F}^\bullet)$ is an inclusion of $A_v$ into $\mathbb{Q}$.*

Now we define a complex $R\,\mathbb{G}_{\mathrm{m}}$, up to quasi-isomorphism, of $Sh(\overline{X}_{\acute{e}t})$ by an exact triangle

$$\mathbb{G}_{\mathrm{m}} \to R\,\mathbb{G}_{\mathrm{m}} \to \mathscr{F}^\bullet[-1] \to,$$

or equivalently an exact triangle

$$R\,\mathbb{G}_{\mathrm{m}} \to \mathscr{F}^\bullet[-1] \xrightarrow{\gamma} \mathbb{G}_{\mathrm{m}}[1] \to .$$

One can choose carefully a morphism $\gamma \in \operatorname{Hom}_{D(Sh(\overline{X}_{\acute{e}t}))}(\mathscr{F}^\bullet[-1], \mathbb{G}_{\mathrm{m}}[1]) \cong \operatorname{Ext}^2_{\overline{X}_{\acute{e}t}}(\mathscr{F}^\bullet, \mathbb{G}_{\mathrm{m}})$ so that $R\,\mathbb{G}_{\mathrm{m}}$ has the cohomologies that we expected. Before computing the group $\operatorname{Ext}^2_{\overline{X}_{\acute{e}t}}(\mathscr{F}^\bullet, \mathbb{G}_{\mathrm{m}})$, we need to derive the following lemmas.

**Lemma 5.2.** *Let $I$ be an inductive system of index that can be refined by an index which is isomorphic to $\mathbb{Z}$, then the canonical morphism $\operatorname{Ext}^s_{\mathcal{C}}(\varinjlim_i A_i, B) \longrightarrow \varprojlim_i \operatorname{Ext}^s_{\mathcal{C}}(A_i, B)$ is an isomorphism if any of the following conditions holds*

(a) *the projective system $\{\operatorname{Ext}^{s-1}_{\mathcal{C}}(A_i, B)\}$ satisfies the Mittag-Leffler condition,*
(b) *$\operatorname{Ext}^{s-1}_{\mathcal{C}}(A_i, B)$ is a finite group for all $i$.*

*Proof.* Because $I$ can be refined by $\mathbb{Z}$, so $\varinjlim_I A_i = \varinjlim_{\mathbb{Z}} A_i$ and $\varprojlim_I \operatorname{Ext}^p_{\mathcal{C}}(A_i, B) = \varprojlim_{\mathbb{Z}} \operatorname{Ext}^p_{\mathcal{C}}(A_i, B)$. Thus, we only need to consider the case that $I \cong \mathbb{Z}$. When $I \cong \mathbb{Z}$, by [12], there is a spectral sequence

$$\varprojlim{}^r_I \operatorname{Ext}^s_{\mathcal{C}}(A_i, B) \Rightarrow \operatorname{Ext}^{r+s}_{\mathcal{C}}(\varinjlim_I A_i, B).$$

This spectral sequence is certainly degenerated as we know that $\varprojlim^r$ is vanishing for any inductive system of abelian groups when $r \geq 2$ (cf. [16] 3.5). Therefore, we get exact sequences

$$0 \to \varprojlim{}^1_I \operatorname{Ext}^{s-1}_{\mathcal{C}}(A_i, B) \to \operatorname{Ext}^s_{\mathcal{C}}(\varinjlim_I A_i, B) \to \varprojlim{}^1_I \operatorname{Ext}^s_{\mathcal{C}}(A_i, B) \to 0,$$

for all $s \geq 0$.

The Mittag-Leffler condition implies $\varprojlim^1 \operatorname{Ext}^{s-1}_{\mathcal{C}}(A_i, B) = 0$ and condition (b) also implies the vanishing of $\varprojlim^1 \operatorname{Ext}^{s-1}_{\mathcal{C}}(A_i, B)$ by [6] 2.3. Thus, the isomorphism is established under both conditions. □

**Lemma 5.3.** (a)

$$\mathrm{Ext}^p_{\overline{X}_{\acute{e}t}}(\mathbb{Q},\mathbb{G}_{\mathrm{m}}) = \begin{cases} 0 & p=2, \\ \widehat{\mathbb{Z}}\bigotimes_{\mathbb{Z}}\mathbb{Q} = \prod'_q \mathbb{Q}_q & p=3, \end{cases}$$

*where* $\prod'_q \mathbb{Q}_q = \mathbb{A}_{\mathbb{Q}}/\mathbb{R}$ *is the restricted product of the* $\mathbb{Q}_p$*'s with respect to the* $\mathbb{Z}_p$*'s.*

(b)

$$\mathrm{Ext}^p_{\overline{X}_{\acute{e}t}}(\alpha_*\mathscr{M}_v,\mathbb{G}_{\mathrm{m}}) = \begin{cases} 0 & p=2, \\ H^3(\overline{X}_{\acute{e}t},\mathbb{G}_{\mathrm{m}}) = \mathbb{Q}/\mathbb{Z} & p=3 \text{ and } v \mid \infty, \\ \varprojlim_n H^3(\overline{X},\mathbb{G}_{\mathrm{m}}) = \widehat{\mathbb{Z}}\bigotimes_{\mathbb{Z}}\mathbb{Q} & p=3 \text{ and } v\nmid\infty. \end{cases}$$

*Proof.* (a) Since $\mathbb{Q} = \varinjlim_n \mathbb{Z}$ and $H^i(\overline{X}_{\acute{e}t},\mathbb{G}_{\mathrm{m}})$ are finite groups for $i=1,2$, we have

$$\mathrm{Ext}^p_{\overline{X}_{\acute{e}t}}(\mathbb{Q},\mathbb{G}_{\mathrm{m}}) = \varprojlim_n H^p(\overline{X}_{\acute{e}t},\mathbb{G}_{\mathrm{m}}).$$

As $H^2(\overline{X}_{\acute{e}t},\mathbb{G}_{\mathrm{m}})=0$, one sees that $\mathrm{Ext}^2_{\overline{X}_{\acute{e}t}}(\mathbb{Q},\mathbb{G}_{\mathrm{m}})=0$. Also, $\varprojlim_n H^3(\overline{X}_{\acute{e}t},\mathbb{G}_{\mathrm{m}}) = \varprojlim_n \mathbb{Q}/\mathbb{Z} = \mathrm{Hom}(\mathbb{Q},\mathbb{Q}/\mathbb{Z}) = \widehat{\mathbb{Z}}\bigotimes_{\mathbb{Z}}\mathbb{Q}$.

(b) As $\mathscr{M}_v = \varinjlim_L \pi'_{L,*}A_v$,

$$\mathrm{Ext}^p_{\overline{X}_{\acute{e}t}}(\alpha_*\mathscr{M}_v,\mathbb{G}_{\mathrm{m}}) = \mathrm{Ext}^p_{\overline{X}_{\acute{e}t}}(\varinjlim_L \pi'_{L,*}A_v,\mathbb{G}_{\mathrm{m}}).$$

When $v$ is an infinite place, $A_v=\mathbb{Z}$, and by the norm theorem(c.f. [9] III.3.9 and [2] Lemma 3.8),

$$\mathrm{Ext}^p_{\overline{X}_{\acute{e}t}}(\pi'_{L,*}\mathbb{Z},\mathbb{G}_{\mathrm{m}}) = H^p(u(L^v)_{\acute{e}t},\mathbb{G}_{\mathrm{m}}),$$

for all $p$. Note that $F$ is totally imaginary and $u(L^v)$ contains all the finite places of $L^v$, so $H^1(u(L^v)_{\acute{e}t},\mathbb{G}_{\mathrm{m}}) = Pic(u(L^v)), H^2(u(L^v)_{\acute{e}t},\mathbb{G}_{\mathrm{m}})=0$, and $H^3(u(L^v)_{\acute{e}t},\mathbb{G}_{\mathrm{m}}) = \mathbb{Q}/\mathbb{Z}$. Since $H^1(u(L^v)_{\acute{e}t},\mathbb{G}_{\mathrm{m}})$ and $H^2(u(L^v)_{\acute{e}t},\mathbb{G}_{\mathrm{m}})$ are both finite groups, by Lemma 5.2, we actually get

$$\mathrm{Ext}^p_{\overline{X}_{\acute{e}t}}(\alpha_*\mathscr{M}_v,\mathbb{G}_{\mathrm{m}}) = \varprojlim_L H^p(u(L^v)_{\acute{e}t},\mathbb{G}_{\mathrm{m}})$$

for $p=2,3$. Clearly, $\mathrm{Ext}^2_{\overline{X}_{\acute{e}t}}(\alpha_*\mathscr{M}_v,\mathbb{G}_{\mathrm{m}}) = H^2(u(L^v)_{\acute{e}t},\mathbb{G}_{\mathrm{m}})=0$. Moreover, for any finite extension $K$ over $L$, the transition map $H^3(u(K^v)_{\acute{e}t},\mathbb{G}_{\mathrm{m}}) \to H^3(u(L^v)_{\acute{e}t},\mathbb{G}_{\mathrm{m}})$ is an isomorphism ([2] Prop.3.2). So, we conclude that $\mathrm{Ext}^3_{\overline{X}_{\acute{e}t}}(\alpha_*\mathscr{M}_v,\mathbb{G}_{\mathrm{m}}) = H^3(\overline{X}_{\acute{e}t},\mathbb{G}_{\mathrm{m}}) = \mathbb{Q}/\mathbb{Z}$.

For the case $v$ is a finite place of $F$, there is a slightly different, as now $A_v = \mathbb{Q} = \varinjlim \frac{1}{n}\mathbb{Z}$ and so $\alpha_*M_v = \varinjlim_{L,n} \pi'_{L,*}\mathbb{Z}$. We claim that

$$\mathrm{Ext}^p_{\overline{X}_{\acute{e}t}}(\varinjlim_{L,n} \pi'_{L,*}\mathbb{Z},\mathbb{G}_{\mathrm{m}}) = \varprojlim_{L,n} H^p(u(L^v)_{\acute{e}t},\mathbb{G}_{\mathrm{m}}),$$

for $p=2,3$. Indeed, because $\mathrm{Ext}^1_{u(L^v)_{\acute{e}t}}(\pi'_{L,*}\mathbb{Z},\mathbb{G}_{\mathrm{m}}) = H^1(u(L^v)_{\acute{e}t},\mathbb{G}_{\mathrm{m}}) = Pic(u(L^v))$ is a finite groups, and $\mathrm{Ext}^2_{u(L^v)_{\acute{e}t}}(\pi'_{L,*}\mathbb{Z},\mathbb{G}_{\mathrm{m}}) = H^2(u(L^v)_{\acute{e}t},\mathbb{G}_{\mathrm{m}}) = 0$, by Lemma 5.2, one gets $\mathrm{Ext}^p_{\overline{X}_{\acute{e}t}}(\alpha_*\mathscr{M}_v,\mathbb{G}_{\mathrm{m}}) = \varprojlim_{L,n} H^p(u(L^v)_{\acute{e}t},\mathbb{G}_{\mathrm{m}})$ for $p=2,3$. Again, by Lemma 5.2, one get $\mathrm{Ext}^p_{\overline{X}_{\acute{e}t}}(M_v,\mathbb{G}_{\mathrm{m}}) = \varprojlim_{L,n} \mathrm{Ext}^p_{u(L^v)}(\mathbb{Z},\mathbb{G}_{\mathrm{m}})$ for $p=2,3$. Hence

$$\mathrm{Ext}^p_{\overline{X}_{\acute{e}t}}(M_v,\mathbb{G}_{\mathrm{m}}) = \begin{cases} 0 & p=2, \\ \varprojlim_n H^3(\overline{X}_{\acute{e}t},\mathbb{G}_{\mathrm{m}}) & p=3. \end{cases}$$

□

**Proposition 5.2.**

$$\mathrm{Ext}^p_{\overline{X}_{\acute{e}t}}(\mathscr{F}^0, \mathbb{G}_\mathrm{m}) = \begin{cases} 0 & p = 2, \\ \prod_{v|\infty} H^3(\overline{X}_{\acute{e}t}, \mathbb{G}_\mathrm{m}) = \prod_{\overline{X}_\infty} \mathbb{Q}/\mathbb{Z} & p = 3. \end{cases}$$

*Proof.* Let $\mathscr{F}^0_v := ker(\alpha_* \mathscr{M}_v \to i_{v,*}A_v)$. Then it arises a long exact sequence

$$\xrightarrow{\partial^{p-1}} \mathrm{Ext}^p_{\overline{X}_{\acute{e}t}}(\alpha_* M_v, \mathbb{G}_\mathrm{m}) \to \mathrm{Ext}^p_{\overline{X}_{\acute{e}t}}(F^0_v, \mathbb{G}_\mathrm{m}) \to \mathrm{Ext}^p_{\overline{X}_{\acute{e}t}}(i_{v,*}A_v, \mathbb{G}_\mathrm{m}) \xrightarrow{\partial^p} \mathrm{Ext}^{p+1}_{\overline{X}_{\acute{e}t}}(\alpha_* M_v, \mathbb{G}_\mathrm{m}) \to .$$

Since $\mathrm{Ext}^2_{\overline{X}_{\acute{e}t}}(\alpha_* \mathscr{M}_v, \mathbb{G}_\mathrm{m}) = 0$ by the previous lemma, and it is easy to see that $\mathrm{Ext}^4_{\overline{X}_{\acute{e}t}}(i_{v,*}A_v, \mathbb{G}_\mathrm{m}) = 0$, we have the following exact sequence,

$$0 \to \mathrm{Ext}^2_{\overline{X}_{\acute{e}t}}(\mathscr{F}^0_v, \mathbb{G}_\mathrm{m}) \to \mathrm{Ext}^3_{\overline{X}_{\acute{e}t}}(i_{v,*}A_v, \mathbb{G}_\mathrm{m}) \to \mathrm{Ext}^3_{\overline{X}_{\acute{e}t}}(\alpha_* \mathscr{M}_v, \mathbb{G}_\mathrm{m}) \to \mathrm{Ext}^3_{\overline{X}_{\acute{e}t}}(\mathscr{F}^0_v, \mathbb{G}_\mathrm{m}) \to 0.$$

By the adjunction of $i_{v,*}$ and $i_v^!$ and Lemma 5.2 (as $H^2_v(\overline{X}_{\acute{e}t}, \mathbb{G}_\mathrm{m}) = 0$),

$$\mathrm{Ext}^3_{\overline{X}_{\acute{e}t}}(i_{v,*}A_v, \mathbb{G}_\mathrm{m}) = \begin{cases} H^3_v(\overline{X}_{\acute{e}t}, \mathbb{G}_\mathrm{m}) = 0 & v \mid \infty, \\ \varprojlim_n H^3_v(\overline{X}_{\acute{e}t}, \mathbb{G}_\mathrm{m}) = \prod'_q \mathbb{Q}_q & v \nmid \infty. \end{cases}$$

Therefore, $\mathrm{Ext}^2_{\overline{X}_{\acute{e}t}}(\mathscr{F}^0_v, \mathbb{G}_\mathrm{m}) = 0$ and $\mathrm{Ext}^3_{\overline{X}_{\acute{e}t}}(\mathscr{F}^0_v, \mathbb{G}_\mathrm{m}) = \mathrm{Ext}^3_{\overline{X}_{\acute{e}t}}(\alpha_* \mathscr{M}_v, \mathbb{G}_\mathrm{m}) = \mathbb{Q}/\mathbb{Z}$ when $v|\infty$.

We claim that the morphism $\mathrm{Ext}^3_{\overline{X}_{\acute{e}t}}(i_{v,*}A_v, \mathbb{G}_\mathrm{m}) \to \mathrm{Ext}^3_{\overline{X}_{\acute{e}t}}(\alpha_* M_v, \mathbb{G}_\mathrm{m})$ is an isomorphism when $v \nmid \infty$. For this, observe that there is a commutative diagram of sheaves

$$\begin{array}{ccc} \alpha_* \mathscr{M}_v = \varinjlim_L \pi'_{L,*} A_v & \longrightarrow & i_{v,*} A_v \\ \uparrow & \nearrow & \\ A_v = \pi'_{F,*} A_v & & \end{array}$$

where $A_v \to i_{v,*}A_v$ is the canonical morphism induced by adjunction of $i_v^*$ and $i_{v,*}$. Applying $\mathrm{Ext}^3(-, \mathbb{G}_\mathrm{m})$ to the above diagram, we get another commutative diagram :

$$\begin{array}{ccc} \mathrm{Ext}^3(\alpha_* \mathscr{M}_v, \mathbb{G}_\mathrm{m}) & \xleftarrow[i]{} & \varprojlim_n H^3_v(\overline{X}_{\acute{e}t}, \mathbb{G}_\mathrm{m}) \\ \downarrow & \swarrow_h & \\ \varprojlim_n H^3(\overline{X}_{\acute{e}t}, \mathbb{G}_\mathrm{m}). & & \end{array}$$

All the groups in the above diagram are isomorphic to $\widehat{\mathbb{Z}} \bigotimes_\mathbb{Z} \mathbb{Q} = \prod'_q \mathbb{Q}_q$. Note the $h$ is an isomorphism as the natural morphism $H^3_v(\overline{X}_{\acute{e}t}, \mathbb{G}_\mathrm{m}) \to H^3(\overline{X}_{\acute{e}t}, \mathbb{G}_\mathrm{m})$ is an isomorphism. Therefore, $i$ is an injection from $\prod'_q \mathbb{Q}_q$ to itself. Since $\mathrm{Hom}_\mathbb{Z}(\mathbb{Q}_p, \mathbb{Q}_q) = 0$ is $p \neq q$ and $\mathbb{Q}_q$ are fields, one sees that $i$ has to be an isomorphism. It follows that $\mathrm{Ext}^2_{\overline{X}_{\acute{e}t}}(\mathscr{F}^0_v, \mathbb{G}_\mathrm{m})$ and $\mathrm{Ext}^3_{\overline{X}_{\acute{e}t}}(\mathscr{F}^0_v, \mathbb{G}_\mathrm{m})$ are vanishing.

Consequently, $\mathrm{Ext}^2_{\overline{X}_{\acute{e}t}}(\mathscr{F}^0, \mathbb{G}_\mathrm{m}) = 0$ and

$$\mathrm{Ext}^3_{\overline{X}_{\acute{e}t}}(\mathscr{F}^0, \mathbb{G}_\mathrm{m}) = \bigoplus_{v \in \overline{X}} \mathrm{Ext}^3_{\overline{X}_{\acute{e}t}}(\mathscr{F}^0_v, \mathbb{G}_\mathrm{m}) = \bigoplus_{v|\infty} \mathbb{Q}/\mathbb{Z}$$

□

**Proposition 5.3.** *There is a canonical isomorphism*

$$\mathrm{Ext}^2_{\overline{X}_{\acute{e}t}}(\mathscr{F}^\bullet, \mathbb{G}_{\mathrm{m}}) \cong \widehat{\mathbb{Z}}.$$

*Proof.* The exact sequence

$$0 \longrightarrow ker \sum \longrightarrow \mathscr{F}^0 \longrightarrow \mathbb{Q} \longrightarrow 0,$$

induces a long exact sequence

$$\mathrm{Ext}^2_{\overline{X}_{\acute{e}t}}(\mathscr{F}^0, \mathbb{G}_{\mathrm{m}})(= 0) \to \mathrm{Ext}^2_{\overline{X}_{\acute{e}t}}(ker \sum, \mathbb{G}_{\mathrm{m}}) \to \mathrm{Ext}^3_{\overline{X}_{\acute{e}t}}(\mathbb{Q}, \mathbb{G}_{\mathrm{m}}) \xrightarrow{\Delta} \bigoplus_{v|\infty} H^3(\overline{X}_{\acute{e}t}, \mathbb{G}_{\mathrm{m}}) \to,$$

where $\Delta$ is a copies of the canonical map $\pi : \mathrm{Hom}(\mathbb{Q}, \mathbb{Q}/\mathbb{Z}) \to \mathrm{Hom}(\mathbb{Z}, \mathbb{Q}/\mathbb{Z})$.

Clearly, $\ker \Delta = \ker \pi = \mathrm{Hom}_{\mathbb{Z}}(\mathbb{Q}/\mathbb{Z}, \mathbb{Q}/\mathbb{Z}) = \widehat{\mathbb{Z}}$. As a consequence $\mathrm{Ext}^2_{\overline{X}_{\acute{e}t}}(ker \sum, \mathbb{G}_{\mathrm{m}}) = ker\Delta = ker(\pi) = \widehat{\mathbb{Z}}$.

□

Let $\gamma$ be the class in $\mathrm{Ext}^2_{\overline{X}_{\acute{e}t}}(\mathscr{F}^\bullet, \mathbb{G}_{\mathrm{m}})$ that corresponds to the generator $\mathbb{1}$ of $\widehat{\mathbb{Z}}$, and $R\,\mathbb{G}_{\mathrm{m}}$ is the complex decided by the exact triangle

$$R\,\mathbb{G}_{\mathrm{m}} \to \mathscr{F}^\bullet[-1] \xrightarrow{\gamma} \mathbb{G}_{\mathrm{m}}[1] \to .$$

To compute the étale cohomology of the complex $R\,\mathbb{G}_{\mathrm{m}}$, we need to determine the coboundary morphisms $\partial^i_U : H^i(U_{\acute{e}t}, \mathscr{F}^\bullet) \xrightarrow{\smile\gamma} H^{i+2}(U_{\acute{e}t}, \mathbb{G}_{\mathrm{m}})$, for any connected étale $U \to \overline{X}$.

**Lemma 5.4.** *For any connected étale $U \to \overline{X}$, let $S := \overline{U} - U$, we have*

(a) *$\partial^0_U$ is the canonical projection*

$$(\bigoplus_{v\in S} A_v)^{\Sigma=0} \to (\bigoplus_{v\in S} A_v/\mathbb{Z})^{\Sigma=0},$$

(b) *$\partial^1_U$ is an isomorphism if $S$ is non-empty, and it is the canonical projection $\mathbb{Q} \to \mathbb{Q}/\mathbb{Z}$ otherwise.*

(c) *$\partial^i_U$ is isomorphic for all $i \geq 2$.*

*In particular, all the $\partial^i$ are surjective for all $i \geq 0$.*

*Proof.* We first deal with the second part of (b). Consider the following commutative diagram

$$\begin{array}{ccccc}
H^1(\overline{X}_{\acute{e}t}, ker \sum) & \times & \mathrm{Ext}^2_{\overline{X}_{\acute{e}t}}(ker \sum, \mathbb{G}_{\mathrm{m}}) & \xrightarrow{\smile} & H^3(\overline{X}_{\acute{e}t}, \mathbb{G}_{\mathrm{m}}) \\
\uparrow\wr & & \downarrow & & \| \\
H^0(\overline{X}_{\acute{e}t}, \mathbb{Q}) & \times & \mathrm{Ext}^3_{\overline{X}_{\acute{e}t}}(\mathbb{Q}, \mathbb{G}_{\mathrm{m}}) & \xrightarrow{\smile} & H^3(\overline{X}_{\acute{e}t}, \mathbb{G}_{\mathrm{m}}) \\
\| & & \| & & \| \\
\varinjlim H^0(\overline{X}_{\acute{e}t}, \frac{1}{n}\mathbb{Z}) & \times & \varprojlim \mathrm{Ext}^3_{\overline{X}_{\acute{e}t}}(\frac{1}{n}\mathbb{Z}, \mathbb{G}_{\mathrm{m}}) & \xrightarrow{\smile} & H^3(\overline{X}_{\acute{e}t}, \mathbb{G}_{\mathrm{m}})
\end{array}$$

It is easy to see that the bottom cup product sends pair $(a, \mathbb{1})$ to $(a \mod \mathbb{Z})$ and thus $a \smile \gamma = a \smile \mathbb{1} = a \mod \mathbb{Z}$. This shows $\partial^1 : H^1(\overline{X}_{\acute{e}t}, \mathscr{F}^\bullet) \to H^3(\overline{X}_{\acute{e}t}, \mathbb{G}_{\mathrm{m}})$ is the canonical projection $\mathbb{Q} \to \mathbb{Q}/\mathbb{Z}$. The general case that $U = \overline{U}$ follows from the commutative diagram

$$\begin{array}{ccc}
H^1(\overline{X}_{\acute{e}t}, \mathscr{F}^\bullet) & \xrightarrow[\text{mod } \mathbb{Z}]{\smile\gamma} & H^3(\overline{X}_{\acute{e}t}, \mathbb{G}_{\mathrm{m}}) \\
\Big\downarrow = & & \Big\downarrow = \\
H^1(\overline{U}_{\acute{e}t}, \mathscr{F}^\bullet) & \xrightarrow{\smile\gamma} & H^3(\overline{U}_{\acute{e}t}, \mathbb{G}_{\mathrm{m}}).
\end{array}$$

For a Zariski open subset $j : U \hookrightarrow \overline{X}$, set $i : \overline{X} - U \hookrightarrow \overline{X}$, we have the following commutative diagram

$$\begin{array}{ccccccccccc}
0 & \longrightarrow & H^0(U_{\acute{e}t}, \mathscr{F}^\bullet) & \longrightarrow & \bigoplus_{v \notin U} H^1_v(\overline{X}_{\acute{e}t}, \mathscr{F}^\bullet) & \xrightarrow{\Sigma} & H^1(\overline{X}_{\acute{e}t}, \mathscr{F}^\bullet) & \longrightarrow & H^1(U_{\acute{e}t}, \mathscr{F}^\bullet) & \longrightarrow & 0 \\
 & & \Big\downarrow \smile\gamma & & \Big\downarrow \smile\gamma & & \Big\downarrow \smile\gamma & & \Big\downarrow \smile\gamma & & \\
0 & \longrightarrow & H^2(U_{\acute{e}t}, \mathbb{G}_{\mathrm{m}}) & \longrightarrow & \bigoplus_{v \notin U} H^3_v(\overline{X}_{\acute{e}t}, \mathbb{G}_{\mathrm{m}}) & \xrightarrow{\Sigma} & H^3(\overline{X}_{\acute{e}t}, \mathbb{G}_{\mathrm{m}}) & \longrightarrow & H^3(U_{\acute{e}t}, \mathbb{G}_{\mathrm{m}}) & \longrightarrow & 0,
\end{array}$$

where the rows are exact sequences induced by the short exact sequence

$$0 \longrightarrow j_! j^* F \longrightarrow F \longrightarrow i_* i^* F \longrightarrow 0, \text{ for any } F \in Sh(\overline{X}).$$

The morphism $H^0(U_{\acute{e}t}, ker \sum) \xrightarrow{\smile\gamma} H^2(U_{\acute{e}t}, \mathbb{G}_{\mathrm{m}})$ is completely determined by $H^1(\overline{X}_{\acute{e}t}, ker \sum) \xrightarrow{\smile\gamma} H^3(\overline{X}_{\acute{e}t}, \mathbb{G}_{\mathrm{m}})$ because $\Sigma$ restricts to an inclusion on $H^1_v(\overline{X}_{\acute{e}t}, \mathscr{F}^\bullet) = A_v$. (This can be seen by put $U$ to be $\overline{X} - v$ and use the fact $H^0((\overline{X} - v)_{\acute{e}t}, ker \sum) = 0$). Therefore, by the commutativity, $\smile \gamma : H^0(U_{\acute{e}t}, \mathscr{F}^\bullet) \to H^2(U_{\acute{e}t}, \mathbb{G}_{\mathrm{m}})$ maps $(\bigoplus_{v \in \overline{U} - U} A_v)^{\Sigma = 0}$ to $(\bigoplus_{v \in \overline{U} - U} A_v/\mathbb{Z})^{\Sigma = 0}$ by taking modular by $\mathbb{Z}$ componentwise.

For general $i : U \to \overline{X}$, set $U^a$ to be the largest open subscheme of $U$ such that for any $v \in U^a$, all the valuations of $K(U)$ above $i(v)$ is contained in $U^a$. We have the commutative diagram

$$\begin{array}{ccc}
H^0(U_{\acute{e}t}, \mathscr{F}^\bullet) & \xrightarrow{\smile\gamma} & H^2(U_{\acute{e}t}, \mathbb{G}_{\mathrm{m}}) \\
\Big\downarrow & & \Big\downarrow \\
H^0(U^a_{\acute{e}t}, \mathscr{F}^\bullet) & \xrightarrow{\smile\gamma} & H^2(U^a_{\acute{e}t}, \mathbb{G}_{\mathrm{m}}) \\
\Big\uparrow & & \Big\uparrow \\
H^0(i(U^a)_{\acute{e}t}, \mathscr{F}^\bullet) & \xrightarrow{\smile\gamma} & H^2(i(U^a)_{\acute{e}t}, \mathbb{G}_{\mathrm{m}})
\end{array}$$

More precisely,

$$\begin{array}{ccc}
\Big(\bigoplus_{w\in\overline{U}-U} A_w\Big)^{\Sigma=0} & \xrightarrow{\smile\gamma} & \Big(\bigoplus_{w\in\overline{U}-U} A_w/\mathbb{Z}\Big)^{\Sigma=0} \\
\big\downarrow & & \big\downarrow \\
\Big(\bigoplus_{w\in\overline{U}-U^a} A_w\Big)^{\Sigma=0} & \xrightarrow{\smile\gamma} & \Big(\bigoplus_{w\in\overline{U}-U^a} A_w/\mathbb{Z}\Big)^{\Sigma=0} \\
\big\uparrow & & \big\uparrow \\
\Big(\bigoplus_{v\in\overline{X}-i(U^a)} A_v\Big)^{\Sigma=0} & \xrightarrow{\smile\gamma} & \Big(\bigoplus_{v\in\overline{X}-i(U^a)} A_v/\mathbb{Z}\Big)^{\Sigma=0}
\end{array}$$

Note that the cup products are canonical on $v$ and $w$, and we have seen above that the bottom cup product is canonical projection. By chasing diagram, one sees that the middle one is also canonical projection, and so is the top one as the injection is canonical. It follows that $\partial_U^0$ is the canonical projection.

To see the first part of part (b), we assume $S$ is non-empty. When $\overline{U}_f \nsubseteq U$, $H^3(U_{ét}, \mathbb{G}_m) = 0$ and $H^1(U_{ét}, \mathscr{F}^\bullet) = 0$ (Prop. 5.1). Therefore, we have $\partial_U^1 = 0$. For the case $\overline{U}_f \subseteq U$, note that $\overline{U}$ is then étale over $\overline{X}$ as all the infinite places are complex. Therefore, we have the following commutation diagram

$$\begin{array}{ccc}
\mathbb{Q} = H^1(\overline{U}_{ét}, \mathscr{F}^\bullet) & \xrightarrow{\pi} & H^1(U_{ét}, \mathscr{F}^\bullet) = \mathbb{Q}/\mathbb{Z} \\
\big\downarrow \partial_{\overline{U}}^1 & & \big\downarrow \partial_U^1 \\
\mathbb{Q}/\mathbb{Z} = H^3(\overline{U}_{ét}, \mathbb{G}_m) & = & H^3(U_{ét}, \mathbb{G}_m) = \mathbb{Q}/\mathbb{Z}.
\end{array}$$

As we have seen, in the very beginning of the proof, that $\partial_{\overline{U}}^1$ and $\pi$ are the canonical projections, $\partial_U^1$ has to be the identity map on $\mathbb{Q}/\mathbb{Z}$ by the commutative.

Part (c) is trivial because $H^i(U_{ét}, \mathscr{F}^\bullet) = H^{i+2}(U_{ét}, \mathbb{G}_m) = 0$ for all $i \geq 2$. □

Now, we are ready to compute the étale cohomology of the complex $R\,\mathbb{G}_m$.

**Proposition 5.4.** *Let $U$ be étale over $\overline{X}$, $S = \overline{U} - U$ and $L = K(U)$, then*

$$H^p(U_{ét}, R\,\mathbb{G}_m) = \begin{cases} \mathcal{O}_{L,S}^\times & p = 0, \\ \mathrm{Pic}(\mathcal{O}_{L,S}) \oplus X_S & p = 1, \\ 0 & p = 2, \\ \mathbb{Z}\,(resp.\,0) & p = 3 \text{ and } S = \emptyset\,(resp.\ S \neq \emptyset), \\ 0 & p \geq 4. \end{cases}$$

*Hence, these cohomologies coincide with those of the $\mathbb{Z}$-dual of $\tau^{\leq 3}\,\mathrm{R}\Gamma_c(U_W, \mathbb{Z})[2]$ (c.f. Prop. 4.5).*

*Proof.* Consider the usual long exact sequence on cohomology induced by the exact triangle $R\,\mathbb{G}_m \to \mathscr{F}^\bullet[-1] \xrightarrow{\gamma} \mathbb{G}_m[1] \to$ and use the result that $\partial^i$ are surjective (Lemma 5.4), we get

$$H^0(U_{ét}, \mathbb{G}_m) \xrightarrow{\sim} H^0(U_{ét}, R\,\mathbb{G}_m),$$

$$0 \to H^1(U_{ét}, \mathbb{G}_m) \to H^1(U_{ét}, R\,\mathbb{G}_m) \to \ker(\partial^0) \to 0 \text{ is exact,}$$

$$H^i(U_{\acute{e}t}, R\,\mathbb{G}_{\mathrm{m}}) = \ker(\partial^{i-1}).$$

Thus, $H^0(U_{\acute{e}t}, R\,\mathbb{G}_{\mathrm{m}}) = \mathcal{O}_{L,S}^{\times}$. Because $\ker\partial^0 = (\oplus_{v\in S}\mathbb{Z})^{\Sigma=0} = X_S$ is $\mathbb{Z}$-free, the middle exact sequence splits and so $H^1(U_{\acute{e}t}, R\,\mathbb{G}_{\mathrm{m}}) = Pic(\mathcal{O}_{L,S}) \oplus X_S$. Also, $H^3(U_{\acute{e}t}, R\,\mathbb{G}_{\mathrm{m}}) = 0$ (resp. $\mathbb{Z}$) when $S \neq \emptyset$ (resp. $S = \emptyset$) and $H^i(U_{\acute{e}t}, R\,\mathbb{G}_{\mathrm{m}}) = 0$ for all $i \geq 4$, by Lemma 5.4 (b) and (c). □

The computations of Prop.5.4 and 4.5 suggest that

**Theorem 5.1.** $\mathrm{R}\Gamma(U_{\acute{e}t}, R\,\mathbb{G}_{\mathrm{m}}) \xrightarrow{\text{qis}} \mathrm{RHom}\left(\tau^{\leq 3}\,\mathrm{R}\Gamma_c(U_W, \mathbb{Z}), \mathbb{Z}\right)[-2]$ *in* $D(\mathbb{Z})$.

In the next section, we shall prove the following more general quasi-isomorphism

$$\mathrm{R}\Gamma(U_{\acute{e}t}, R\,\mathbb{G}_{\mathrm{m}}) \xrightarrow{\text{qis}} \mathrm{RHom}\left(\tau^{\leq 3}\,\mathrm{R}\Gamma_c(U_W, \mathbb{Z}), \mathbb{Z}\right)[-2],$$

in $D(\mathbb{Z}[G])$, when $S$ is stable under the action of $G$. Note that when $G$ is the trivial group, this is just Theorem 5.1

5.2. **The Duality Theorem.** Throughout this section, we assume that $F$ is a *totally imaginary* number field, $X = \operatorname{Spec}\mathcal{O}_F$, $U = \operatorname{Spec}\mathcal{O}_{L,S}$ be any connected étale neighborhood of $X$. We also assume that $K$ is a subfield of $L$ such that the extension $L/K$ is Galois and $G := \mathrm{Gal}(L/K)$. We require that $S$ contains all the archimedean places and is stable under the action of $G$.

Recall that $\mathscr{F}^\bullet(U) \simeq \ker\sum(U)[0]$ in $D(\mathbb{Z}[G])$. To simplify the notation, in the following, we denote by $\mathscr{F}^\bullet(U)$ the complex $\ker\sum(U)[0]$.

**Theorem 5.2.** *Let* $U = \operatorname{Spec}\mathcal{O}_{L,S}$ *and* $A \to B$ *represents* $\mathrm{RHom}(\tau^{\leq 3}\,\mathrm{R}\Gamma_c(U_W, \mathbb{Z}), \mathbb{Z})[-2]$ *in* $\mathrm{Ext}^2_G(X_S \oplus \mathrm{Pic}(U), U_S)$. *Then* $\mathrm{R}\Gamma(U_{\acute{e}t}, \mathbb{G}_{\mathrm{m}})$ *can be represented by*

$$A \to B \to \mathscr{F}^\bullet(U).$$

*Moreover, there is an exact triangle in* $D(\mathbb{Z}[G])$,

$$\mathscr{F}^\bullet(U)[-2] \to \mathrm{R}\Gamma(U_{\acute{e}t}, \mathbb{G}_{\mathrm{m}}) \to \mathrm{RHom}(\tau^{\leq 3}\,\mathrm{R}\Gamma_c(U_W, \mathbb{Z}), \mathbb{Z})[-2] \to . \tag{12}$$

*This exact triangle is the same as the exact sequence of complexes*

$$0 \longrightarrow \mathscr{F}^\bullet(U)[-2] \longrightarrow (A \to B \to \mathscr{F}^\bullet(U)) \longrightarrow (A \to B) \longrightarrow 0.$$

*Proof.*

$$\begin{array}{ccccccc}
\uparrow & & \uparrow & & \uparrow & & \\
\mathrm{RHom}(\tau^{\leq 3}\,\mathrm{R}\Gamma_c(U_W, \mathbb{Z}), \mathbb{Z})[-2] & \longrightarrow & \mathrm{RHom}(\mathrm{R}\Gamma_c(U_{\acute{e}t}, \mathbb{Z}), \mathbb{Z})[-2] & \longrightarrow & I_0[0] & \longrightarrow & \\
\vdots\!\uparrow & & \uparrow & & \uparrow & & \\
\mathrm{R}\Gamma(U_{\acute{e}t}, \mathbb{G}_{\mathrm{m}}) & \longrightarrow & \mathrm{R}\Gamma_c(U_{\acute{e}t}, \mathbb{Z})^D[-3] & \xrightarrow{\alpha} & I_0[0] \bigoplus I_2[-2] & \longrightarrow & \\
\vdots\!\uparrow & & \uparrow & & \uparrow & & \\
\mathscr{F}^\bullet(U)[-2] & \longrightarrow & X_S \bigotimes_{\mathbb{Z}} \mathbb{Q}[-2] & \longrightarrow & I_2[-2] & \longrightarrow & , \\
 & & \| & & & & \\
 & & \Big(\mathrm{RHom}(\mathrm{R}\Gamma_c(U_{\acute{e}t}, \mathbb{Z}), \mathbb{Q})[-3]\Big) & & & &
\end{array}$$

where $I_0 = \widehat{U}_S/U_S$ and $I_2 = \bigoplus_{v\in S_\infty} \mathbb{Q}/\mathbb{Z}$ are both $\mathbb{Z}[G]$-injective modules. All the rows are known exact triangles. The top is applying $\mathrm{RHom}_{\mathbb{Z}}(-,\mathbb{Z})$ to the exact triangle

$$\mathrm{R}\Gamma_c(U_{\acute{e}t},\mathbb{Z}) \to \mathrm{R}\Gamma_c(U_W,\mathbb{Z}) \to \mathrm{Hom}_{\mathbb{Z}}(U_S,\mathbb{Q})[-2] \to,$$

the middle is (5) induced by the Artin-Verdier Duality Theorem, and the bottom is followed from the definition of $\mathscr{F}^\bullet(U)$ and is indeed a short exact sequence. The middle column is obtained by applying

$$\mathrm{Hom}_{\mathbb{Z}}(-,\mathbb{Z}) \to \mathrm{Hom}_{\mathbb{Z}}(-,\mathbb{Q}) \to \mathrm{Hom}_{\mathbb{Z}}(-,\mathbb{Q}/\mathbb{Z}) \to$$

to the complex $\mathrm{R}\Gamma_c(U_{\acute{e}t},\mathbb{Z})$.

The top right hand and bottom right hand squares of the above diagram are commutative because they are commutative on $H^0$ and $H^2$ respectively (see Lemma 4.2). Thus, by property of derived categories, there exist morphisms in $D(\mathbb{Z}[G])$ so that the first column is an exact triangle and the above 3 by 3 diagram is semi-commutative.

Suppose $\mathrm{RHom}(\tau^{\leq 3}\,\mathrm{R}\Gamma_c(U_W,\mathbb{Z}),\mathbb{Z})[-2]$ is in the same class of $A \xrightarrow{f} B$ as an element of the 2-extension group $\mathrm{Ext}^2_G(X_S \oplus \mathrm{Pic}(U), U_S)$. Then $\mathrm{RHom}(\tau^{\leq 3}\,\mathrm{R}\Gamma_c(U_W,\mathbb{Z}),\mathbb{Z})[-2]$ is in the class of image, says $A' \xrightarrow{f} B$, of $A \to B$ in $\mathrm{Ext}^2_G(X_S \oplus \mathrm{Pic}(U), \widehat{U}_S)$, and $\mathrm{R}\Gamma_c(U_{\acute{e}t},\mathbb{Z})^D[-3]$ is in the same class of $A' \to B \to X_S \bigotimes_{\mathbb{Z}} \mathbb{Q}$. In fact, one may choose $A'$ to be $(A \oplus \widehat{U}_S)/U_S$ and set $f'(\overline{a},x) = f(a)$. Note that $\mathrm{R}\Gamma(U_{\acute{e}t},\mathbb{G}_{\mathrm{m}})$ is determined by $\alpha$, more precisely by $H^0(\alpha)$ and $H^2(\alpha)$. It's easy to see that the exact triangle associated to the following exact sequence of complexes is in fact isomorphic to the middle row of the above diagram.

$$\begin{array}{ccccc}
0 & & 0 & & 0 \\
\downarrow & & \downarrow & & \downarrow \\
A & \xrightarrow{f} & B & \longrightarrow & \mathscr{F}(U) \\
\downarrow & & \downarrow & & \downarrow \\
A' & \xrightarrow{f'} & B & \longrightarrow & X_S \bigotimes_{\mathbb{Z}} \mathbb{Q} \\
\downarrow & & \downarrow & & \downarrow \\
I_0 & \longrightarrow & 0 & \longrightarrow & I_2 \\
\downarrow & & \downarrow & & \downarrow \\
0 & & 0 & & 0
\end{array}$$

We conclude that $\mathrm{R}\Gamma(U_{\acute{e}t},\mathbb{G}_{\mathrm{m}})$ can be represented by $A \to B \to \mathscr{F}(U)$ and it induced the exact triangle

$$\mathscr{F}^\bullet(U)[-2] \to \mathrm{R}\Gamma(U_{\acute{e}t},\mathbb{G}_{\mathrm{m}}) \to \mathrm{RHom}_{\mathbb{Z}}(\tau^{\leq 3}\,\mathrm{R}\Gamma_c(U_W,\mathbb{Z}),\mathbb{Z})[-2] \to .$$

$\square$

**Remark.** *By definition of the complex $R\,\mathbb{G}_{\mathrm{m}}$, $\mathrm{R}\Gamma(U_{\acute{e}t},R\,\mathbb{G}_{\mathrm{m}})$ satisfies the following exact triangle*

$$\mathscr{F}^\bullet(U)[-2] \to \mathrm{R}\Gamma(U_{\acute{e}t},\mathbb{G}_{\mathrm{m}}) \to \mathrm{R}\Gamma(U_{\acute{e}t},R\,\mathbb{G}_{\mathrm{m}}) \to . \tag{13}$$

*One expects that exact triangles (12) and (13) are isomorphic, which would imply the quasi-isomorphism*

$$\mathrm{RHom}(\tau^{\leq 3}\,\mathrm{R}\Gamma_c(U_W,\mathbb{Z}),\mathbb{Z})[-2] \simeq \mathrm{R}\Gamma(U_{\acute{e}t}, R\,\mathbb{G}_{\mathrm{m}}).$$

**Theorem 5.3.** *For any $U = \operatorname{Spec}\mathcal{O}_{L,S}$, one has a quasi-isomorphism*

$$\mathrm{RHom}(\tau^{\leq 3}\,\mathrm{R}\Gamma_c(U_W,\mathbb{Z}),\mathbb{Z})[-2] \simeq \mathrm{R}\Gamma(U_{\acute{e}t}, R\,\mathbb{G}_{\mathrm{m}}),$$

*in $D(\mathbb{Z}[G])$.*

*Proof.* We first show that

$$\mathrm{Hom}_{D(\mathbb{Z}[G])}(\mathscr{F}^\bullet(U)[-2], \mathrm{R}\Gamma(U_{\acute{e}t},\mathbb{G}_{\mathrm{m}})) \simeq \mathrm{Hom}_G(\mathscr{F}^\bullet(U), H^2(U_{\acute{e}t},\mathbb{G}_{\mathrm{m}})) \tag{14}$$

Observe that there is an exact sequence

$$0 \to X_{S_f}\otimes_{\mathbb{Z}}\mathbb{Q} \to \mathscr{F}^\bullet(U) \to \bigoplus_{S_\infty}\mathbb{Z} \to 0.$$

Since $X_{S_f}\otimes_{\mathbb{Z}}\mathbb{Q}$ and $\bigoplus_{S_\infty}\mathbb{Z}\,(\cong \mathbb{Z}[\mathrm{Gal}(F/\mathbb{Q})])$ are $G$-cohomologically trivial, so is $\mathscr{F}^\bullet(U)$. By lemma 4.2(b), there is an exact sequence

$$0 \longrightarrow \mathrm{Ext}^1_{\mathbb{Z}[G]}(\mathscr{F}^\bullet(U), H^1(U_{\acute{e}t},\mathbb{G}_{\mathrm{m}})) \longrightarrow \mathrm{Hom}_{D(\mathbb{Z}[G])}(\mathscr{F}^\bullet(U)[-2], \mathrm{R}\Gamma(U_{\acute{e}t},\mathbb{G}_{\mathrm{m}}))$$

$$\longrightarrow \mathrm{Hom}_G(\mathscr{F}^\bullet(U), H^2(U_{\acute{e}t},\mathbb{G}_{\mathrm{m}})) \longrightarrow 0$$

Note that $\mathrm{Ext}^1_G(\mathscr{F}^\bullet(U), H^1(U_{\acute{e}t},\mathbb{G}_{\mathrm{m}})) = \mathrm{Ext}^1_G(X_{S_f}\otimes_{\mathbb{Z}}\mathbb{Q}, \mathrm{Pic}(U)) = 0$. Hence the isomorphism (14).

By Lemma 5.4 and by Theorem 5.2, we see that the induced homomorphisms on $H^2$,

$$\mathscr{F}^\bullet(U) \to H^2(U_{\acute{e}t},\mathbb{G}_{\mathrm{m}}),$$

for (12) and (13) are both the canonical projections. Thus, exact triangles (12) and (13) have to be isomorphic, in virtue of isomorphism (14). Consequently,

$$\mathrm{RHom}_{\mathbb{Z}}(\tau^{\leq 3}\,\mathrm{R}\Gamma_c(U_W,\mathbb{Z}),\mathbb{Z})[-2] \simeq \mathrm{R}\Gamma(U_{\acute{e}t}, R\,\mathbb{G}_{\mathrm{m}}),$$

in $D(\mathbb{Z}[G])$.

□

**Corollary 5.2.** *Suppose that $\mathrm{Pic}(\mathrm{Spec}(\mathcal{O}_{F,S})) = 0$ and $S$ is stable under the action of G, then we have the following quasi-isomorphism*

$$\mathrm{R}\Gamma((\mathrm{Spec}(\mathcal{O}_{F,S}))_{\acute{e}t}, R\,\mathbb{G}_{\mathrm{m}}) \simeq \Psi_S,$$

*where $\Psi_S$ is the canonical Tate sequence associated to $S$.*

*Proof.* This follows from Theorem 5.3 and Theorem 4.3. □

**Remark.** (a) *Theorem 5.3 allows us to recover the Weil-étale cohomology groups $H^p_c(U_W,\mathbb{Z})$ by the hyper étale cohomology of the complex $R\,\mathbb{G}_{\mathrm{m}}$. Also, one might find corresponding Weil-étale cohomology axioms in terms of $H^p(U_{\acute{e}t}, R\,\mathbb{G}_{\mathrm{m}})$.*

(b) *The reason that it is hard to generalize Lichtenbaum's prototype to higher dimensional arithmetic schemes $\mathscr{X}$ is that there are no Weil groups for higher dimensional fields. However, Theorem 5.3 shows us a probability to generalize Lichtenbaum's prototype, because we do not use Weil groups when defining the complex $R\,\mathbb{G}_{\mathrm{m}}$. One direct thought is that, for any $n$-dimensional arithmetic scheme $\mathscr{X}$, one may define a complex $R\mathbb{Z}(n)$ in $\mathrm{Ext}^2_{\overline{\mathscr{X}}}(\mathbb{Z}(n), \mathscr{F}^\bullet(n))$, so that the $\mathbb{Z}$-dual of $\mathrm{R\Gamma}(\overline{\mathscr{X}}_{et}, R\mathbb{Z}(n))$ defines certain Weil-étale cohomology theory, where $\mathscr{F}^\bullet(n)$ is a complex of étale sheaves that depends on $n$.*

(c) *Theorem 5.3 suggests that for any étale $U \to \overline{X}$ there could be a perfect pairing*

$$\mathrm{RHom}_{U_{\acute{e}t}}(\mathscr{F}, R\,\mathbb{G}_{\mathrm{m}}) \times \tau^{\leq 3}\,\mathrm{R\Gamma}_c(U_W, \mathscr{F}) \to \mathbb{Z},$$

*for a certain class of $U_{\acute{e}t}$-sheaves $\mathscr{F}$.*

## References

[1] E. Artin and J. Tate, *Class Field Theory*, Benjamin, New York, 1967.

[2] Mel Bienenfeld, *An étale cohomology duality theorem for number fileds with a real embedding*, Trans. A.M.S. Vol. 303(Sep. 1987), No.1, pp. 71-96,

[3] D. Burns and M. Flach,*On Galois structure invariants associated to Tate motives*, Amer. J. Math. 120 (1998), 1343-1397.

[4] M. Flach, *Cohomology of topological groups with applications to the Weil group*, Compositio Math., Vol. 144, no. 3, 2008.

[5] M. Flach and B. Morin, *On the Weil-étale topos of arithmetic schemes*, to appear.

[6] C. Jensen, *Les Foncteurs Dérivés de* $\varprojlim$ *et leurs Applications en Théorie des Modules,* Lecture Notes in Math., 254, Springer, Heidelberg.

[7] F. Knudsen and D. Mumford, *The projectivity of the moduli space of stable curves. I : Preliminaries on "det" and "Div"*, Math. Scand. 39(1976), 19-55.

[8] S. Lichtenbaum, *The Weil-étale topology for number rings*, Ann. of Math. Vol. 170 (2009), No. 2, 657-683.

[9] J. S. Milne, *Aritmetic Duality Theorems*, BookSurge, LLC, 2nd Edition, 2006.

[10] B. Morin, Sur le topos Weil-étale dun corps de nombres, Ph.D. thesis, l'université Bordeaux I, 2008.
`Http://grenet.drimm.u-bordeaux1.fr/pdf/2008/MORIN_BAPTISTE_2008.pdf`

[11] B. Morin, *On the Weil-étale cohomology of number fields*, Trans. Amer. Math. Soc. 363 (2011), 4877-4927.

[12] J-E. Roos, *Sur les foncteurs dérivés de* $\varprojlim$. Applications. Compte Rendu Acad. Sci. Paris 252, 3702-3704. Erratum ibid. 254 p1722.

[13] J. Tate, *Number theoretic background* , Proceedings of Symposia in Pure Mathematics, Vol. 33 (1979), part 2, 3-26.

[14] J.Tate, *Les conjectures de Stark sur les fonctions L d'Artin en s=0*, Birkhäuser, 1984.

[15] J. L. Verdier, *Des Catégories Dérivées des Categories Abéliennes*, Asterisque 239 (1996).

[16] C. A. Weibel, *An Introduction to Homological Algebra,* Cambridge Studies in Advanced Math., vol. 38, Cambridge Univ. Press, 1994.

[17] T. Zink, *étale cohomology and duality in number fields*, Haberland, *Galois cohomology*, Berlin, 1978, Appendix 2.

California Institute of Technology,, 1200 E. California Blvd. Pasadena, CA 91125.
*E-mail address*: `ycchiu@caltech.edu`